\documentclass[10pt,letterpaper]{amsart}
\usepackage{color, enumerate}
\usepackage{graphicx}
\usepackage{amssymb}
\usepackage{bbm, hyperref}
\usepackage{amsmath}

\newtheorem{thm}{Theorem}

\newtheorem{cor}[thm]{Corollary}
\newtheorem{prop}[thm]{Proposition}
\newtheorem{lemma}[thm]{Lemma}
\theoremstyle{definition}
\newtheorem{example}[thm]{Example}
\newtheorem{remark}[thm]{Remark}
\newtheorem{definition}[thm]{Definition}
\newtheorem{question}[thm]{Question}
\newcommand{\norm}[1]{\left\vert #1 \right \vert}

\numberwithin{thm}{subsection}
\numberwithin{equation}{section}

\newcommand{\R}{\mathbb R}
\newcommand{\C}{\mathbb C}
\newcommand{\Z}{\mathbb Z}
\newcommand{\ii}{\mathbbm{i}}
\newcommand{\Heis}{\mathbb H}
\newcommand{\Sp}{\mathbb S}

\newcommand{\st}{:}

\def\Barint_#1{\mathchoice
          {\mathop{\vrule width 6pt height 3 pt depth -2.5pt
                  \kern -8pt \intop}\nolimits_{#1}}%
          {\mathop{\vrule width 5pt height 3 pt depth -2.6pt
                  \kern -6pt \intop}\nolimits_{#1}}%
          {\mathop{\vrule width 5pt height 3 pt depth -2.6pt
                  \kern -6pt \intop}\nolimits_{#1}}%
          {\mathop{\vrule width 5pt height 3 pt depth -2.6pt
                  \kern -6pt \intop}\nolimits_{#1}}}

\title{Heisenberg quasiregular ellipticity}
\author[K. F\"assler]{Katrin F\"assler}
\address{Department of Mathematics and Statistics \\ University of Jyv\"askyl\"a \\ P.O.Box 35 (MaD), FI-40014 University of Jyv\"askyl\"a, Finland\\ \emph{Current address: }Department of Mathematics\\ University of Fribourg\\Chemin du Mus\'{e}e 23,
CH-1700 Fribourg, Switzerland}
\email{katrin.faessler@unifr.ch}
\author[A. Lukyanenko]{Anton Lukyanenko}
\address{Department of Mathematics \\ University of Michigan \\ 530 Church Street, Ann Arbor, MI 48109}
\email{anton@lukyenenko.net}
\author[J. T. Tyson]{Jeremy T. Tyson}
\address{Department of Mathematics \\ University of Illinois \\ 1409 West Green St. \\ Urbana, IL, 61801}
\email{tyson@illinois.edu}
\date{\today}
\thanks{K.F.\ was supported by the Academy of Finland through the grant 285159 `Sub-Riemannian manifolds from a quasiconformal viewpoint'. A.L.\ was supported by NSF RTG grant  DMS-1045119. J.T.T.\ was supported by Simons Collaboration Grant 353627 `Geometric analysis in sub-Riemannian and metric spaces' and NSF grants DMS-1201875 `Geometric mapping theory in sub-Riemannian and metric spaces'. and DMS-1600650 `Mappings and measures in sub-Riemannian and metric spaces'.}
\keywords{Quasiregular mapping, contact manifold, sub-Riemannian manifold, $3$-sphere, link complement, Hopf link, isoperimetric inequality, Sobolev--Poincar\'e inequality, capacity, nonlinear potential theory, morphism property}

\begin{document}

\begin{abstract}
Following the Euclidean results of Varopoulos and Pankka--Rajala, we
provide a necessary topological condition for a sub-Riemannian
3-manifold $M$ to admit  a nonconstant quasiregular
mapping from
the sub-Riemannian Heisenberg group $\mathbb{H}$. As an application, we show that a link complement
$\Sp^3\backslash L$ has a sub-Riemannian metric admitting such a mapping  only if $L$ is empty, the unknot or Hopf link. In the converse direction, if $L$ is empty, a specific unknot or Hopf link, we construct a quasiregular mapping from $\mathbb{H}$ to $\Sp^3\backslash L$.

The main result is obtained by translating a growth condition on $\pi_1(M)$ into the existence of a supersolution to the $4$-harmonic equation, and relies on recent advances in the study of analysis and potential theory on metric spaces.
\end{abstract}

\maketitle

\section{Introduction}

\setcounter{subsection}{0}
\setcounter{thm}{0}

Given two topological manifolds $X$ and $Y$, it is often quite difficult to decide whether there exists a covering map from $f: X\rightarrow Y$.
The only obvious obstruction is that the universal covers of $X$ and $Y$ should be homeomorphic.
Furthermore, two manifolds with the same universal cover may have substantially different geometries.
For example, $\R^2$ covers both the torus with abelian fundamental group $\Z^2$, and the punctured torus with fundamental group the free group $\mathbb F_2$.

The covering problem becomes more tractable if we impose geometric
restrictions on the allowed covering maps.  We will ask for $f$ to
be \emph{quasiregular}, imposing a quasiconformal-type restriction
on metric distortion, but also allowing the mapping to be a
{branched} covering map {onto its image}; or (for compactness
purposes) a constant mapping. See Section \ref{s:prelim} for a precise
definition of quasiregularity. Examples of quasiregular mappings
include isometric embeddings, conformal and quasiconformal
home\nobreak{}omorphisms, and branched holomorphic mappings.

In this paper, we ask which sub-Riemannian 3-manifolds admit quasiregular mappings from the Heisenberg group, that is, which such manifolds are {\it Heisenberg quasiregularly elliptic}. We begin by reviewing the history of the quasiregular ellipticity question.

\subsection{Euclidean quasiregular ellipticity}
A Riemannian $n$-manifold $M$ is said to be \emph{(Euclidean)
quasiregularly elliptic} if there exists a nonconstant
quasiregular mapping from $\R^n$ to $M$  (see \S \ref{sec:DefEllipticity} for a more nuanced discussion of the definition).
The classification of
quasiregularly elliptic manifolds goes back to the Picard Theorem
in complex analysis, and has had a profound impact on the
development of geometric mapping theory. In particular, Gromov
devotes a chapter of his celebrated `Green Book', \cite[Chapter
6]{MR2307192} to the interplay between isoperimetric inequalities,
quasiregular ellipticity and the geometry of groups.

The classical Picard Theorem states that every nonconstant entire function misses at most one point. By Sto\"ilow
factorization \cite{MR0082545,MR0344463},
this is equivalent to saying that $\R^2 \setminus P$ is
quasiregularly elliptic if and only if $P$ contains at most one
point. While the holomorphic interpretation does not persist in
higher dimensions, Rickman showed in \cite{MR583633} for all $n\geq 3$
that if $\R^n \setminus P$ is quasiregularly elliptic, then $P$
contains at most finitely many points. Rickman also provided a converse in three dimensions
\cite{MR781587}. It was only in a remarkable recent paper of
Drasin and Pankka \cite{MR3372169} that Rickman's construction was extended
to all dimensions. Alternate PDE proofs of the Rickman--Picard
Theorem were provided by Lewis and Eremenko--Lewis
\cite{MR1195483, MR1139803}.

Varopoulos \cite[pp.\
146--147]{MR1218884} proved that if a closed Riemannian
$n$-manifold $M$ is quasiregularly elliptic, then the fundamental
group $\pi_1(M)$  is virtually nilpotent and indeed has growth
rate at most $n$. Pankka--Rajala \cite{MR2832708} extended
Varopoulos' theorem to open manifolds, and provided a result in
the spirit of Picard's Theorem: a link complement $\Sp^3\backslash
L$ is \emph{quasiregularly elliptic for some choice of Riemannian
metric} if and only if $L$ is empty, the unknot or the Hopf link
(cf.\ Gromov \cite[Examples 6.12]{MR2307192}). The present paper
generalizes the results of Pankka--Rajala and thus contributes to the study of quasiregular ellipticity in a sub-Riemannian setting.

In the Riemannian setting, the theory is naturally more advanced.
 A full classification of
closed quasiregularly elliptic 3-manifolds was provided by
Jormakka in \cite{MR973719}. Holopainen and Rickman
\cite{MR1191342} extended the Rickman-Picard theorem to more
general Riemannian targets,
and they studied quasiregular mappings {between} Riemannian
manifolds in \cite{MR1667818}. Bonk and Heinonen provided an
obstruction to quasireguar ellipticity  for closed manifolds in
terms of cohomological dimension in \cite{MR1846030}, which allows
to prove nonellipticity in some cases where the fundamental group
is too small to apply Varopoulos' theorem.

\subsection{Heisenberg quasiregular ellipticity}
In this paper, we leave the Riemannian framework and study a
quasiregular ellipticity in the
{sub}-Riemannian setting (\S \ref{ss:SRmfds}). The simplest
homogeneous space admitting a
non-Euclidean sub-Riemannian metric is the Heisenberg group
$\Heis$, a nilpotent group of step 2 and topological dimension 3.
The Heisenberg group shares much of the structure of Euclidean
space, including a one-parameter family of metric dilations, and
therefore serves as a natural model space for sub-Riemannian
geometry and source space for quasiregular mappings. We say that a
sub-Riemannian manifold is \emph{Heisenberg quasiregularly
elliptic} if it admits a non-constant quasiregular map from
$\Heis$.

The study of Heisenberg
quasiregular ellipticity is in the early stages of development.
A Rickman--Picard theorem for the Heisenberg group (and more generally for H-type Carnot groups) was provided by Heinonen--Holopainen in \cite{MR1630785} (see also \cite{MR2252688}): if $\Heis \setminus P$ is Heisenberg-quasiregularly elliptic, then $P$ contains at most finitely many points. However, it remains an open problem to show that a quasiregular map from $\Heis$ to itself can miss even a single point.

The study of quasiregular mappings to more general sub-Riemannian targets was recently initiated in \cite{FLP,GNW,GL,GW}.
Following the above results of Varopoulos and Pankka--Rajala, we prove

\begin{thm}
\label{thm:linkOurs} Let $L\subset \Sp^3$ be a link in the
three-sphere and let $\Heis$ denote the first Heisenberg group
equipped with its standard sub-Riemannian structure.
\begin{itemize}
\item If there
exists an equiregular sub-Riemannian metric $g$ on
$\Sp^3\backslash L$ admitting a nonconstant quasiregular mapping
$f: \Heis \rightarrow (\Sp^3\backslash L, g)$ then $L$
is empty, an unknot, or a Hopf link.
\item Conversely, there exist a smooth unknot $S$ and a smooth Hopf link $H$, and for $L \in \{\emptyset, S, H\}$, there exist equiregular sub-Riemannian metrics $g_\emptyset$, $g_S$, and $g_H$ in $\Sp^3\backslash L$ and nonconstant quasiregular maps $f: \Heis \rightarrow (\Sp^3\backslash L, g_L)$.
\end{itemize}
\end{thm}

For the second part, we provide in Section
\ref{ss:ex} new explicit examples of mappings from the Heisenberg
group onto the $3$-sphere and (specific) unknot and Hopf link
complements. The first statement of the theorem is a consequence of the
following more general Varopoulos-type result, which we prove in
Section \ref{s:main}.

\begin{thm}
\label{thm:main} Let $M$ be an equiregular sub-Riemannian
$3$-manifold. If there exists a nonconstant quasiregular mapping
$f: \Heis \rightarrow M$, then the growth rate of $\pi_1(M)$ is at
most $4$.
\end{thm}

Here and throughout this paper, we say that a group $G$ has {\it growth rate larger than $d$} if there exists a finite set $S$ in $G$ and a constant $c>0$ so that the cardinality of any ball $B(R)$ is at least $cR^d$ for all positive integers $R$, where $B(R)$ denotes the ball of radius $R$ about the identity element in the word metric on the subgroup $\langle S\rangle$ generated by $S$. Equivalently, $G$ has {\it growth rate at most $d$} if for every finitely generated subgroup $\Gamma$ of $G$ and for every finite set $S$ with $\Gamma = \langle S\rangle$, there exists a constant $C>0$ so that $B(R)$ has cardinality at most $CR^d$ for all positive integers $R$. See \cite{MR2832708} or \cite[\S5B]{MR2307192} for more details.

\begin{remark}
We expect that a statement as in Theorem \ref{thm:main} holds true also in higher dimensions, that is, for quasiregular maps from $\Heis^n$ to an equiregular sub-Riemannian $(2n+1)$-manifold $M$, with an analogous proof. For simplicity, we restrict our discussion to $3$-manifolds.
\end{remark}

Assuming Theorem \ref{thm:main} we now indicate how to derive Theorem \ref{thm:linkOurs}.

\begin{proof}[Proof of Theorem \ref{thm:linkOurs}]
If $L$ is not one of the links listed above, then $\pi_1(M)$ contains a free group of rank at least $2$, and thus it has exponential growth, that is, in particular it has growth rate larger than $4$. See the references in \cite{MR2832708}. By Theorem  \ref{thm:main}, manifolds with this property cannot admit nonconstant quasiregular mappings from the Heisenberg group.

The examples in Section \ref{ss:ex} establish the positive implication in all the remaining cases, that is, if $L$ is empty, a specific unknot, or a specific Hopf link.
\end{proof}

\subsection{Outline of the proof of Theorem \ref{thm:main}}

The proof of Theorem \ref{thm:main} will be given in Section \ref{s:main}. Here we provide a brief outline, following the corresponding subsections of Section \ref{s:main}.
\begin{enumerate}
\item[\ref{sec:topology}.]\label{st:1} Starting with the assumption that $M$ has a fundamental group with a finitely generated subgroup $\Gamma$ with growth rate larger than 4, we define a ``(relatively) compact core'' $M'\subset M$ and a lift $\widetilde M''$ of $M'$ to $\widetilde M$, satisfying $\widetilde M''/\Gamma = M'$. We define a distance on the closure of ${M'}$, lift it to $\widetilde{M}''$ and show that the resulting space is quasi-isometric to $\Gamma$.
\item[\ref{ss:isop}.] Using the fact that $\widetilde{M}''$ is quasi-isometric to $\Gamma$ and that $\Gamma$ has growth rate larger than $4$, we show that $\widetilde{M}''$ (or rather, a net $Y$ on $\widetilde{M}''$) satisfies a `rough' $d$-dimensional isoperimetric inequality for some $d>4$.
\item[\ref{sec:localisoperimetric}.]\label{st:3} We use the local geometry of $\widetilde{M}''$ to prove a weak $(\frac{4}{3},1)$-Poincar\'{e} (or Sobolev-Poincar\'{e}) inequality and a weak relative $4$-dimensional isoperimetric inequality (for balls {of fixed size} centred in $Y$). This requires a careful study of $\widetilde{M}''$ at and near its boundary.
\item[\ref{sec:globalisoperimetric}.]\label{st:4} Combining the rough and the relative isoperimetric inequality, we deduce that $\widetilde{M}''$ also fulfills a `smooth' $d$-dimensional isoperimetric inequality. We formulate this implication in an axiomatic way, so that it applies also in a more abstract setting.
\item[\ref{ss:cap}] We next study the $4$-capacity in $\widetilde M$ of a ball in $\widetilde M''$. Fixing an admissible function $u$ for the capacity, we restrict $u$ to $\widetilde M''$ and use
 a sub-Riemannian coarea formula to relate the horizontal gradient of $u\vert_{\widetilde M''}$ to the perimeter of its level sets. Coupled with the isoperimetric inequality established above, we obtain a uniform positive lower bound for the $L^4$-norm of the horizontal gradient of $u$. That is, we show that $\widetilde M$ is $4$-hyperbolic. (In fact, we show a stronger version of $4$-hyperbolicity of $\widetilde M''$ and combine this with the fact that the inclusion $\widetilde M'' \hookrightarrow \widetilde M$ is bi-Lipschitz to obtain hyperbolicity of $\widetilde M$.)
\item[\ref{sec:potentialtheory}.] We conclude from the
$4$-hyperbolicity of $\widetilde M$ the existence of a positive
nonconstant supersolution to the $4$-harmonic equation, and
the existence of a Green's function for the sub-elliptic $4$-Laplacian at every
point of $\widetilde M$. \item[\ref{s:morphism}.] We show that quasiregular mappings have a morphism
property: the pullback of a supersolution to the $4$-harmonic equation is
a supersolution to a nonlinear operator of type $4$.
\item[\ref{sec:victory}.] Lastly, we suppose that $f: \Heis
\rightarrow M$ is a nonconstant quasiregular map. We lift it to a
nonconstant quasiregular map $\widetilde{f}: \Heis \rightarrow
\widetilde M$ and pull back a nonconstant supersolution to the
$4$-harmonic equation, contradicting the $4$-parabolicity of the
Heisenberg group.
\end{enumerate}

While the preceding outline of the proof is largely the same as in the Riemannian case \cite{MR2832708}, the sub-Riemannian geometry enters the picture in a non-trivial way in most of the steps described above. For instance, the contact structure prevents us from constructing a  {double} of $\overline{M'}$ as in \cite{MR2832708}. We address this issue by carefully analyzing intrinsic balls in $M'$ and we state properties of the intrinsic distance, which we believe to be of independent interest. Furthermore, unlike in \cite{MR2832708}, we cannot apply directly the work by Kanai \cite{MR792983}, which has been formulated for Riemannian manifolds with Ricci curvature bounds. We take the opportunity to translate his argument into a more general axiomatic framework, which applies to our setting. The proof works in metric measure spaces with a mild condition on the volumes of balls. Throughout the individual steps of the proof, we also combine results that have been recently developed in various areas of sub-Riemannian geometry, such as classifications of uniform and Sobolev-Poincar\'{e} domains, notions and properties of horizontal perimeter, and others.

\begin{remark}
One could bypass the discussion of the morphism property by proving a capacity inequality for arbitrary condensers in $\widetilde{M}$. This is the approach employed by Varopoulos, see \cite[Chapter X]{MR1218884}. We expect that a similar argument works in the present setting.
However, since the notion of $\mathcal{A}$-harmonic functions has
classically strong connections with questions of quasiregular
ellipticity and is of independent interest for further
developments, we decided to follow a different route in the
present paper.
\end{remark}

\subsection*{Structure of the paper}
This paper is organized as follows. In Section \ref{s:context-and-examples} we exhibit examples of nonconstant quasiregular mappings from the Heisenberg group onto the $3$-sphere and onto the complement of the unknot and Hopf link. Section \ref{s:prelim} contains background information about quasiregular mappings of sub-Riemannian contact manifolds. Section \ref{s:main}
is the heart of the paper. Here we prove Theorem \ref{thm:main} following the outline previously indicated.
We have relegated to an appendix (Appendix \ref{appendix}) several
basic properties of the calculus of horizontal derivatives.

\subsection*{Acknowledgements} We would like to thank Chang-Yu Guo and Pekka Pankka for discussions related to the subject of this paper. Research for this paper was completed during visits of various subsets of the authors to the University of Bern, the University of Jyv\"askyl\"a and the University of Illinois. The hospitality of all of these institutions is appreciated.

\section{Examples of Heisenberg quasiregularly elliptic spaces}\label{s:context-and-examples}
In this section we describe the Heisenberg group and some spaces
that admit quasiregular mappings from it.

\subsection{Sub-Riemannian manifolds}\label{ss:SRmfds}

Recall that a \emph{sub-Riemannian manifold} is a triple $(M,HM,g_M)$, where $M$ is a connected smooth manifold, $HM\subset TM$ is a smooth bracket-generating distribution, and $g_M$ is a
metric on $HM$. An absolutely continuous curve $\gamma$ in $M$ is \emph{horizontal} if $\dot \gamma$ is almost always in the horizontal distribution $HM$. By Chow's Theorem, any two points of $M$ are connected by a horizontal curve, and one defines the Carnot-Carath\'eodory distance between two points $p, q\in M$ as the infimum of $g_M$-lengths of the horizontal curves joining $p$ to $q$. A sub-Riemannian manifold is furthermore \emph{equiregular} if the distribution $HM$ and its iterated brackets are, in fact, subbundles of $TM$ of constant dimension.

In this paper, we restrict our attention to equiregular
sub-Riemannian manifolds $M$ of dimension 3, and assume that $HM
\neq TM$. It is easy to see that the
bracket-generating condition is
then equivalent to $HM$ being a \emph{contact distribution}. That is,
locally $HM$ is the kernel of a smooth \emph{contact form} $\alpha$
satisfying $\alpha \wedge d \alpha \neq 0$. In particular, the
Darboux Theorem states that locally $(M, HM)$ is
contactomorphic to
the Heisenberg group with its standard contact structure. Note
that this
contactomorphism need not send the metric $g_M$ to the
Heisenberg metric $g_\Heis$.

\begin{example}[Heisenberg group]
In \emph{exponential coordinates of the first kind}, the Heisenberg group $\Heis$ is
given by $\R^3$ with group structure
$$(x,y,t)*(x',y',t') = (x+x', y+y', t+t'-2xy'+2yx').$$
The standard contact form on $\Heis$ is given by
$$\alpha_{\Heis} = dt + 2(x\,dy - y\,dx).$$
Notice that $\alpha_{\Heis}$ is invariant under left
translations. The horizontal distribution $H\Heis$ on $\Heis$ is given by $\ker \alpha_{\Heis}$, and is spanned by the left-invariant vector fields
\begin{equation}\label{eq:horiz_vfd}
X= \partial_x + 2y \, \partial_t\quad\text{and}\quad Y = \partial_y -2x \, \partial_t.
\end{equation}
The sub-Riemannian path metric $d_{\Heis}$ on $\Heis$ is induced
by the inner product $g_{\Heis}$ defined by the line element
$ds_{\Heis}^2 = dx^2+dy^2$ on $H\Heis$.
\end{example}

\subsection{Heisenberg quasiregularly elliptic spaces}\label{ss:ex}
\label{sec:newexamples}

Our main theorem shows that not every
equiregular sub-Riemannian 3-manifold is
Heisenberg quasiregularly elliptic.
We now describe several Heisenberg quasiregularly elliptic spaces.
The first of these is well-known, while we believe that the
remaining constructions are new. We will leave the formal
definition of quasiregularity for the next section (Definition
\ref{d:QR_def}), as all the mappings we mention -- apart from
Example \ref{ex:sphere2} -- are covering mappings that are either
locally isometric or conformal.
 For the
moment, it is sufficient to think of a quasiregular map
$f:(\Heis,H\Heis,g_{\Heis})\to (M,HM,g_M)$ as a continuous
branched cover with the property that $f_{\ast}(H\Heis)\subset HM$ and so that there exists a constant $K$ such that for
almost every $p\in \Heis$, one has
\begin{displaymath}
\sup g_M(f_{\ast} v, f_{\ast} v)\leq K^2 \inf g_M(f_{\ast} v,
f_{\ast} v),
\end{displaymath}
where the supremum and the infimum are taken over all horizontal
vectors $v\in H_p \Heis$ with $g_{\Heis}(v,v)=1$.

\begin{example}[Sub-Riemannian $3$-sphere]
\label{ex:sphere1}
Consider the $3$-sphere $\Sp^3$, viewed as the unit sphere in $\C^2$. The tangent space to $\Sp^3$ at a point $p$ is the real orthogonal complement of the normal vector $\vec{n}(p) = p$. This tangent space is not invariant under multiplication by the imaginary unit $\ii$. The subbundle of the tangent bundle which is invariant under multiplication by $\ii$ coincides with the kernel of a contact form and defines a sub-Riemannian structure. Explicitly, let $\alpha_{\Sp^3}$ be the contact form given by
$$
\alpha_{\Sp^3} = \overline w_1 \, dw_1 - w_1 \, d\overline{w_1} + \overline{w_2} \, dw_2 - w_2 \, d\overline{w_2}
$$
where $w = (w_1,w_2)$ denote coordinates in $\C^2$. The standard
sub-Riemannian metric on $\Sp^3$ is given by the restriction of
the Euclidean inner product to $\ker \alpha_{\Sp^3}$.

The inverse stereographic projection
$$
\iota(x,y,t) = \left( \frac{2y-2x\ii}{1+x^2 + y^2 - \ii t}, \frac{1-x^2 - y^2 + \ii t}{1+x^2 + y^2 - \ii t}\right)
$$
provides a bijection between $\Heis$ and $\Sp^3\setminus\{(0,-1)\}$, and is furthermore well-known to be both contactomorphic and conformal, see for instance \cite[p.\ 315]{KoRe1985} or \cite[Section 3.3]{CDPT}. It follows that the sub-Riemannian $\Sp^3$ and the punctured sub-Riemannian $\Sp^3$ are Heisenberg quasiregularly elliptic.
\end{example}

\begin{example}[Lens spaces]
Let $\mathbb L$ be a lens space equipped with its standard contact structure and sub-Riemannian metric arising from its representation as a quotient of $\Sp^3$.
(See \cite[Section 3.1]{FLP} for details on the sub-Riemannian structure on lens spaces). We assume that $\mathbb L \ne \Sp^3$.
Composing the embedding $\iota: \Heis \hookrightarrow \Sp^3$ with the quotient projection $\pi: \Sp^3\rightarrow \mathbb L$ gives a conformal map from $\Heis$ onto $\mathbb L$.

More generally, if $\Gamma$ is any group of isometries of the sub-Riemannian $\Sp^3$ such that $\Sp^3/\Gamma$ is a smooth manifold and $\pi:\Sp^3 \to \Sp^3/\Gamma$ denotes the quotient map, then the composition $\pi \circ \iota$ is a quasiregular mapping from $\Heis$ to $\mathbb \Sp^3/\Gamma$ with its standard contact structure and sub-Riemannian metric. See \cite{MR2051684} and
\cite{MR2782691} for examples of finite isometry groups of $\Sp^3$ arising in the study of proper holomorphic mappings between balls and CR representation theory.
\end{example}

\begin{example}[Unknot complement]
\label{ex:unknot} Let $L_1\subset \Sp^3$ be an unknot, whose
complement $M_1:=\Sp^3\backslash L_1$ is diffeomorphic to
$\Sp^1\times \R^2$. We claim that $M_1$ has a sub-Riemannian
metric admitting a surjective quasiregular map from $\Heis$.

Consider first the quotient $M'_1 := \Heis/\langle (0,0,1)\rangle$ of the Heisenberg
group by integer translations along the $t$-axis. Since vertical translations
are isometries of $\Heis$, the sub-Riemannian metric
$g_{\Heis}$ projects to a well-defined
sub-Riemannian metric
 on $M'_1$, with the projection map $\pi: \Heis \rightarrow
M_1'$ a surjective local isometry.

Note now that $M_1$ and $M'_1$ are both diffeomorphic to $\Sp^1\times \R^2$, and let $f: M'_1 \to M_1$ be a diffeomorphism. Give $M_1$ the contact structure and sub-Riemannian metric induced by this diffeomorphism. Then the map $f\circ \pi: \Heis \rightarrow M_1$ is a surjective local isometry, as desired.

For an explicit example, let $L'$ denote the $t$-axis in $\Heis$,
 and
consider the mapping $h(x,y,t) = (\cos(2\pi t) e^x, \sin(2\pi t)
e^x, y)$ from $\Heis$ to itself. This mapping commutes with
integer translations along the $t$-axis, and so induces a
sub-Riemannian metric $g$ on the $\Heis \setminus L'$. Then
$h:(\Heis,g_{\Heis}) \to (\Heis\setminus L',g)$ is a quasiregular
surjection.
\end{example}

\begin{example}[Hopf link complement]
Let $L_2\subset \Sp^3$ be the Hopf link, with
$M_2:=\Sp^3\backslash L_2$ diffeomorphic to $\R\times \Sp^1\times
\Sp^1$. We claim that $M_2$  has a sub-Riemannian metric admitting
a surjective quasiregular map from $\Heis$.

Note first that the integer group $\Z^2$ acts on $\Heis$ by the isometries
$$(0,b,c)*(x,y,t) = (x, y+b, t+c+2bx).$$
The quotient space $M_2'=\Heis/\Z^2$ then inherits a sub-Riemannian structure from $\Heis$, and the projection map $\pi: \Heis \rightarrow M_2'$ is a surjective local isometry.

Note that $M_2$ and $M_2'$ are both diffeomorphic to
$\mathbb{R}\times \Sp^1\times \Sp^1$, and let $f: M_2' \rightarrow
M_2$ be a diffeomorphism. Give $M_2$ the contact structure and
sub-Riemannian metric induced by this diffeomorphism. Then the map
$f\circ \pi: \Heis \rightarrow M_2$ is a surjective local
isometry, as desired.

An explicit example is easy to construct as in Example \ref{ex:unknot}, taking $L$ to be the union of the equators in $\Sp^3\subset \C^2$.
\end{example}

\begin{example}[Surjection to the sub-Riemannian $3$-sphere]
\label{ex:sphere2}

We conclude this section by providing a \emph{surjective} quasiregular mapping from the Heisenberg group to the sub-Riemannian 3-sphere.

Fix an integer $a>1$ and let $f_a': \C^2 \rightarrow \C^2$ be the continuous extension of the map
$$
f'_a:(r_1 e^{\ii \theta_1}, r_2 e^{\ii \theta_2}) \mapsto (r_1 e^{a \ii \theta_1}, r_2 e^{a\ii \theta_2}),
$$
defined on $\Omega_0 := \{(r_1 e^{\ii \theta_1}, r_2 e^{\ii
\theta_2})\,:\,r_1 \ne 0,r_2\ne 0\}$.

The \emph{multi-twist map} $f_a: \Sp^3 \rightarrow \Sp^3$ is given
by restricting $f'_a$ to the sphere. It is clear that $f_a$ is a
branched covering map, with branching along the standard
Hopf link, and it was shown in \cite{FLP} that $f_a$ is
quasiregular according to the so-called metric definition. The
equivalence of various definitions of quasiregularity (metric,
geometric, and analytic) on equiregular sub-Riemannian manifolds
has been established in \cite{GNW,GL}. To keep our discussion
self-contained, we verify here directly that the map appearing in
the following lemma is quasiregular in the sense of Definition
\ref{d:QR_def}.

\begin{lemma}\label{claim:sphereSurjection} Let $\iota$ be the conformal embedding of $\Heis$ into $\Sp^3$ as in Example \ref{ex:sphere1}.
Then the map $f:=f_a \circ \iota: \Heis \to \Sp^3$ is a surjective
quasiregular map. Here $\Heis$ and $\Sp^3$ are endowed with their
standard contact structures and sub-Riemannian metrics.
\end{lemma}

\begin{proof}
Since we have chosen $a$ to be an integer larger than $1$, the map
$f$ is surjective from $\Heis$ to $\Sp^3$. Recall that $\iota$ is
a diffeomorphism and $f_a$ is smooth on $\Omega_0 =\{(r_1 e^{\ii
\theta_1}, r_2 e^{\ii \theta_2})\,:\,r_1 \ne 0,r_2\ne 0\}$, so $f$
is smooth on $\iota^{-1}(\Omega_0)$. By a short computation as in
\cite[Section 3.2]{FLP}, one sees that
\begin{displaymath}
g_{\Sp^3}(v,v) \leq g_{\Sp^3}((f_a)_{\ast} v, (f_a)_{\ast} v) \leq
a^2 g_{\Sp^3}(v,v),\quad v\in H_p \Sp^3, \, p\in \Omega_0.
\end{displaymath}
Combined with the fact that $\iota$ is conformal, it follows that
$f= f_a \circ \iota$ fulfills the distortion estimate required for
quasiregularity on the set $\iota^{-1}(\Omega_0)$. As
$\Heis\setminus \iota^{-1}(\Omega_0)$ is the union of the $t$-axis
with a unit circle in the $xy$-plane, it is negligible and we know
that the distortion estimate holds almost everywhere on $\Heis$ as
required.

The remaining property in Definition \ref{d:QR_def} to verify  is the existence of weak horizontal derivatives of $f$ in
$L^4_{loc}$. Since $f$ is smooth outside  the $t$-axis and a
planar circle in $\Heis$, it follows that it is absolutely
continuous along almost every fiber in a fibration given by a
horizontal left invariant vector field. The corresponding
directional derivatives exist pointwise almost everywhere and are
weak derivatives. (See, for instance, \cite[Theorem
2.2]{MR1404088} and Remark \ref{r:ACL_char}.) It remains to
establish their local integrability. For this purpose, it is
useful to introduce the complex operators $Z=\frac{1}{2}(X-\ii Y)$
and $\bar Z=\frac{1}{2}(X+\ii Y)$ on $\Heis$, and to use
coordinates $w=(w_1,w_2)$ on $\Sp^3$. A direct computation gives
\begin{equation}\label{eq:pushforward_vfd}
\iota_{\ast} Z = - \ii \frac{(1+w_2)^2}{1+\bar w_2} W\quad
\text{and}\quad \iota_{\ast} \bar Z = \ii \frac{(1+\bar
w_2)^2}{1+w_2}\bar W,
\end{equation}
where $W= \bar w_2 \partial_{w_1} - \bar w_1 \partial_{w_2}$ and
$\bar W = w_2 \partial_{\bar w_1} - w_1 \partial_{\bar w_2}$; cf.\
\cite[p.\ 320]{KoRe1985}. Writing $w_j = r_j e^{\ii \theta_j}$ for
$j \in \{1,2\}$ and $w'=(w_1',w_2') = f_a(w)$, we find
$\partial_{w_1} w_1' = e^{\ii
(a-1)\theta_1}\left(\frac{1+a}{2}\right)$ and $\partial_{\bar w_1}
w_1' = e^{\ii (a+1)\theta_1}\left(\frac{1-a}{2}\right)$ on
$\Omega_0$. Analogous formulae hold for $w_2'$. It follows that
$Z(h\circ f_a \circ \iota)$ and $\bar Z (h\circ f_a \circ \iota)$
are in $L^{\infty}_{loc}$ for an arbitrary smooth function
$h:\Sp^3 \to \mathbb{R}$. This shows that $f$ has weak horizontal
derivatives in $L^4_{loc}$ and concludes the proof.
%
%
%
\end{proof}
\end{example}

\subsection{Notions of quasiregular ellipticity}
\label{sec:DefEllipticity}
In this section we discuss some subtleties in the definition of quasiregular ellipticity, and provide a few more examples.

We start with Euclidean quasiregular ellipticity. A Riemannian $n$-manifold $(M, g_M)$ is \emph{quasiregularly elliptic} if there exists a non-constant quasiregular mapping $f: \R^n \rightarrow M$.
A differentiable $n$-manifold $M$ (without specifying a Riemannian metric)  is \emph{quasiregularly elliptic as  a manifold} if it supports some Riemannian metric $g'$  such that $(M, g')$ is quasiregularly elliptic.

\begin{example}
The Rickman--Picard theorem states that any quasiregular map $f:\Sp^3 \to \Sp^3$ of the Euclidean sphere misses at most
finitely many points. Thus, if $L$ is a non-empty link, $\Sp^3 \setminus L$ is not quasiregularly elliptic with the induced standard metric from $\Sp^3$.
However, if $L$ is either a smooth unknot or a smooth Hopf link, then $\Sp^3\setminus L$ is quasiregularly elliptic {as a manifold} \cite{MR2832708}.
\end{example}

Analogous definitions apply for sub-Riemannian 3-manifolds with $\Heis$ as the source space: one can consider Heisenberg quasiregular ellipticity of the sub-Riemannian manifold $(M, HM, g_M)$, of the contact manifold $(M, HM)$ with {some} choice of Riemannian tensor on $HM$, or of the manifold $M$ with {some} choice of bundle and Riemannian tensor. In the examples above, we have shown that the  sub-Riemannian
$3$-sphere and lens spaces are quasiregularly elliptic with their standard sub-Riemannian metric, while the unknot complement and Hopf
link complement are quasiregularly elliptic {as manifolds}.

\begin{question}
Is there an unknot complement or a Hopf link complement that is Heisenberg quasiregularly elliptic when endowed with the contact structure induced by the standard contact form $\alpha_{\Sp^3}$?
\end{question}

We conclude with four more examples of quasiregularly elliptic manifolds, leaving the constructions of quasiregular mappings to the reader.

\begin{example}\label{ex:half_space}
The half-space $\Heis_{x^+} = \{ (x,y,t)\in \Heis \st x>0\}$ is
Heisenberg quasiregularly elliptic {as a contact manifold}.
One can construct an explicit mapping using polarized coordinates.
(For the definition of polarized coordinates, see for instance
\cite[\S 2.1]{CDPT}).
\end{example}

\begin{example}
Any domain $\Omega \subset \Heis$ diffeomorphic to $\Heis$ is Heisenberg quasiregularly elliptic {as a contact manifold}. This follows from the uniqueness of tight contact structures on $\R^3$, see Eliashberg \cite{MR1162559}.
\end{example}

In the next two examples, our source space is $\Heis_{x^+}$.

\begin{example}
Equip $M=\Heis\setminus \{(0,0,t)\st t\in \R\}$ with the standard contact structure. Then there is a Riemannian tensor on $HM$
for which there exists a quasiregular map $f:\Heis_{x^+}\rightarrow \Heis\setminus \{(0,0,t)\st t\in \R\}$. One
can construct an explicit map using polarized coordinates. The map $f$ is not defined on all of $\Heis$ and the contactomorphism from Example \ref{ex:half_space} distorts the standard metric.
Hence the existence of $f$ does not imply
that the unknot complement in $\Sp^3$ with the induced standard
contact structure is Heisenberg quasiregularly elliptic as a contact manifold.
\end{example}

\begin{example}
Let $M=\Sp^3\setminus L$, where $L$ is the union of equators
$w_1=0$ and $w_2=0$ in $\Sp^3$. Give $M$ the standard contact
structure $HM$. Then there is a Riemannian tensor on $HM$ for which there exists a quasiregular map $f:\Heis_{x^+}\rightarrow M$. One can construct an explicit map using a contactomorphism from $\Heis$ to the rototranslation group
(as given for instance in \cite{FKD}). Note that this does not
imply that the Hopf link complement in $\Sp^3$ is Heisenberg
quasiregularly elliptic as a contact manifold with the standard
contact structure, as $f$ is not defined on all of $\Heis$.
\end{example}

\section{Quasiregular mappings: preliminaries}\label{s:prelim}

Quasiregular mappings have first been studied in Euclidean space
as a generalization of complex analytic functions. They also arise
naturally as non-injective counterparts for quasiconformal maps.
We refer to \cite{MR994644,MR1238941} for in-depth introductions
to the subject. The definitions generalize to Riemannian
manifolds, see for instance \cite{MR2492501}. In the
sub-Riemannian setting, quasiregular mappings were first
investigated in the Heisenberg group and other Carnot groups
\cite{MR1630785}, \cite{Da}. Many properties of quasiregular
mappings in Euclidean spaces carry over to the Carnot group
setting; for an example we mention the results in
\cite{MR2252688}. The theory for more general equiregular
sub-Riemannian manifolds has been initiated in \cite{FLP} and
further developed in a recent series of papers \cite{GNW},
\cite{GL},  \cite{GW}, \cite{MR3456891}. Versions of
quasiregularity on metric spaces of locally bounded geometry were
discussed in \cite{MR2275345}.

Even though quasiregular mappings have been studied already in
greater generality, we decided to include in this section a
self-contained discussion, which focuses on the specific setting
of this paper.
Restricting our attention to contact $3$-manifolds allows us to
exploit properties of the Heisenberg group using contact geometry.
Moreover, oriented sub-Riemannian
contact $3$-manifolds can be endowed with a CR structure
\cite{FGV}. While we do not make use of this structure in our
proofs, the reader  interested in CR geometry may read our
results as a continuation of the research on quasiconformal maps
in CR $3$-manifolds
\cite{MR1055843,MR1404088,MR1340848,MR1246889}.


Sections \ref{ss:contactSR} and \ref{ss:defQR} contain definitions
related to sub-Riemannian contact manifolds and quasiregular
mappings, respectively. Section \ref{ss:equiv_char} is devoted to
the interplay between contact geometry and quasiregular mappings,
and we provide auxiliary results that allow us to make use of the
rich theory in the Heisenberg group.

\subsection{Definition of contact forms and measures adapted to a metric}\label{ss:contactSR}
 As discussed in Section \ref{ss:SRmfds}, a $3$-manifold $M$ with a subbundle $HM \subsetneq TM$ is an
 equiregular sub-Riemannian manifold precisely if it is a contact manifold.
 We use this fact throughout the paper, and we alternate our perspective between sub-Riemannian and contact geometry.
We do not assume that the horizontal distribution of the manifolds under consideration is the kernel of an a priori given contact form.
 Instead, we will now describe how to choose a specific such form $\alpha_M$ canonically associated with the sub-Riemannian metric $g_M$ on $(M,HM)$. The related volume form $\alpha_M\wedge \mathrm{d}\alpha_M$ is useful since a meaningful geometric study of quasiregular mappings requires a canonical choice of measure in order to define notions such as Jacobian and distributional Laplacians. In Riemannian manifolds, the choice is the Riemannian volume form. A natural generalization to the sub-Riemannian setting is provided by the
\emph{Popp volume}: a smooth volume form canonically associated to an equiregular sub-Riemannian manifold $M$ and the metric $g_M$ thereon.

\begin{definition}\label{d:Popp} Let $(M,HM,g_M)$ be an equiregular $3$-manifold with a co-orientable horizontal distribution.
Let $\alpha_M$ be the contact form uniquely determined by the conditions that $\mathrm{ker}\,\alpha_M = HM$ and that
$\mathrm{d}{\alpha_M}|_{HM}$ coincides with the volume form induced by $g_M$ on $HM$.
The \emph{Popp volume} $\mathrm{vol}_M$ is given by $\alpha_M\wedge \mathrm{d}\alpha_M$. The associated measure $\mu_M$ is called the \emph{Popp measure}.
 \end{definition}

Equivalently, one can choose a local orthonormal frame $\{e_1,e_2\}$
on $(HM ,g_M)$ and let $e_3$ be the Reeb vector field determined by $\alpha_M$.  If  $\{\nu_1,\nu_2,\nu_3\}$ denotes the dual orthonormal basis to  $\{e_1,e_2,e_3\}$, then $\nu_1\wedge \nu_2 \wedge\nu_3$ agrees with  $\mathrm{vol}_M$, independently of the choice we made for the orthonormal frame.
One can also see $\mathrm{vol}_M$ as the Riemannian volume associated to the extended Riemannian metric obtained from $g_M$ by declaring $\{e_1,e_2,e_3\}$ orthonormal. For a more thorough discussion of Popp measures on contact manifolds, the reader may consult for instance  \cite{MR3263162,MR3207127}.


\begin{remark}\label{r:popp_global_orientable}

In order to keep the presentation in the Definition \ref{d:Popp} simple, we have assumed that the subbundle $HM$ is co-orientable, that is, given by the kernel of a globally defined contact form. If this is not the case, we cannot choose a global orientation for $HM$ and the form $\alpha_M$ is defined only up to a sign. However, while we cannot globally promote $g_M$ to a Riemannian metric, the Popp volume $\mathrm{vol}_M$ is still well-defined. See also \cite[Remark 9]{MR2502528}. In fact, Popp measures can be introduced much more generally on arbitrary equiregular sub-Riemannian manifolds \cite[10.6]{MR1867362}.
\end{remark}

For an equiregular sub-Riemannian $3$-manifold
the Popp measure equals a constant multiple of spherical
$4$-dimensional Hausdorff measure with respect to the sub-Riemannian distance; see \cite[Theorem 4]{MR2875644}.
In the case
of $\Heis$ with the standard sub-Riemannian metric
$g_{\Heis}$, the measure $\mu_{\Heis}$ coincides  with the Haar measure on
$\Heis$, which is the $3$-dimensional Lebesgue measure (all up
to a possible multiplicative factor).

\subsection{Definition of horizontal derivatives and quasiregularity}\label{ss:defQR}
We begin our formal discussion of quasiregularity by introducing
certain classes of functions. Let $U$ be an open set in an
equiregular sub-Riemannian $3$-manifold $N$ with an orthonormal
frame $\{e_1,e_2\}$ of the subbundle $HN$. The \emph{horizontal
Sobolev space} $HW^{1,4}(U)$ is defined as the space of functions
$u\in L^4(U)$ whose distributional derivatives in direction $e_1$
and $e_2$ exist and belong to $L^4(U)$, and the \emph{local
horizontal Sobolev space} $HW_{loc}^{1,4}(U)$ is defined
accordingly. See for instance \cite[Section 2.2]{CLDO} for the
precise definitions.

We also consider the regularity of mappings that take values in a
manifold. However, we do not need the full structure of a Sobolev
space in this setting, so we confine ourselves to the following
definition. It is stated in terms of the (divergence-free) vector
fields $X$ and $Y$ given in \eqref{eq:horiz_vfd}.

\begin{definition}
We say that a continuous map $f:\Heis \to M$ has
\emph{$L^4_{loc}$ weak horizontal derivatives} if for any smooth
function $h:M \to \mathbb{R}$ and any open set $U \Subset
\Heis$, there is a function $g\in L^4_{loc}(\Heis)$
so that
\begin{displaymath}
 \int_U X\varphi \cdot (h\circ f) \; \mathrm{d}\mu_{\Heis} = - \int_U \varphi \cdot g \; \mathrm{d}\mu_{\Heis},\quad \text{for all }\varphi\in \mathcal{C}^{\infty}_0(U),
\end{displaymath}
and analogously for the vector field $Y$. We write $g = X(h \circ
f)$.
\end{definition}

A few remarks concerning this definition are in order.

\begin{remark}
Our definition of weak horizontal derivatives essentially agrees with the one given by Tang in \cite[\S2]{MR1404088}
for maps between smooth strongly pseudoconvex CR $3$-manifolds,
and in case the target manifold is the Heisenberg group, it
matches the standard definition employed in connect with the horizontal Sobolev space
as for instance used in \cite{Da}.
\end{remark}

\begin{remark}\label{r:ACL_char}
According to \cite[Theorem 2.2]{MR1404088}, a continuous map
$f:\Heis \to M$ has  {$L^4_{loc}$ weak horizontal derivatives} if
and only if it is $\mathrm{ACL}$, i.e., $f$ is absolutely
continuous along almost every fiber in the fibration given by $X$
and $Y$, and, moreover, it has $L^4_{loc}$ horizontal derivatives
in these directions. In this case, the weak and pointwise
horizontal derivatives coincide almost everywhere. Analogous
definitions  and statements apply for other integrability
exponents $1\leq p <\infty$.
\end{remark}

We employ the horizontal derivatives to introduce a formal horizontal differential.
Given local coordinates $(x_1,x_2,x_3)$ on $M$, we can define for
every continuous map $f:\Heis \to M$ with $L^4_{loc}$ weak
horizontal derivatives the following notions:
\begin{equation}\label{eq:def_Xf_Yf}
 Xf =\sum_{i=1}^3X (x_i \circ f)\partial_{x_i} \quad \text{and}
 \quad
Yf =\sum_{i=1}^3Y (x_i \circ f)\partial_{x_i}
\end{equation}
The vector fields $Xf$ and $Yf$ are well defined almost
everywhere. Indeed, given charts $\Phi$ and $\Psi$, we write
\begin{displaymath}
X(\Phi_i \circ f) = X((\Phi_i \circ \Psi^{-1})\circ (\Psi \circ f))
\end{displaymath}
and apply the chain rule from Proposition \ref{p:chain} in the
Appendix with $h=\Phi_i $. This yields
\begin{align*}
Xf = \sum_{i=1}^3 X(\Phi_i \circ f) \partial_{\Phi_i}&=
\sum_{i=1}^3 \left( \sum_{j=1}^3 \frac{\partial (\Phi_i \circ
\Psi^{-1})}{\partial \Psi_j}(\Psi(f(\cdot))) X(\Psi_j \circ f)
\right)
\partial_{\Phi_i}\\
&= \sum_{j=1}^3 X(\Psi_j \circ f) \sum_{i=1}^3 \frac{\partial
(\Phi_i \circ \Psi^{-1})}{\partial \Psi_j}\partial_{\Phi_i}\\
&= \sum_{j=1}^3 X(\Psi_j \circ f) \partial_{\Psi_j},
\end{align*}
and analogously for $X$ replaced by $Y$. Thus we can interpret $Xf(p)$ and $Yf(p)$ as elements in $T_p M$ for almost every $p\in M$ and formulate the following definition, which generalizes \cite[Definition 1.1]{Da}.

\begin{definition}
We say that a continuous map $f:\Heis \to M$ with $L^4_{loc}$ weak horizontal derivatives is \emph{weakly contact} if $Xf$ and $Yf$ lie in $HM$ almost everywhere.
\end{definition}

\begin{definition}
Assume that $f:\Heis \to M$ is weakly contact. At almost every $p\in \Heis$, the \emph{formal horizontal differential} of $f$,
\begin{displaymath}
D_H f(p): H_p \Heis \to H_p M,
\end{displaymath}
is defined by
\begin{equation}\label{eq:def_diff}
D_H f(p) (X_p) = Xf(p)\quad\text{and} \quad  D_H f(p) (Y_p) = Y f(p),
\end{equation}
extended to the entire horizontal plane $H_p \Heis$ by linearity.
\end{definition}

To formulate distortion conditions for a weakly contact map
$f:\Heis \to M$, it is convenient to work with the quantities
$\|D_H f(p)\|$ and $\ell[D_H f(p)]$. These are standard notations
in quasiconformal analysis. For a linear map $A:(H_p
\Heis,g_{\Heis}) \to (H_{q} M,g_{M})$, we set
\begin{displaymath}
 \|A\|:= \sup_{g_{\Heis}(v,v)=1}\sqrt{g_M(Av,Av)}\quad \text{and} \quad \ell[A]:= \inf_{g_{\Heis}(v,v)=1}\sqrt{g_M(Av,Av)}.
\end{displaymath}

We are now prepared to state an analytic definition for quasiregularity.

\begin{definition}\label{d:QR_def}
Let $M$ be a smooth orientable $3$-manifold endowed with an equiregular
distribution $HM$ and a sub-Riemannian metric $g_M$. We call a map
$f:\Heis \to M$ \emph{quasiregular} if
\begin{itemize}
\item $f$ is continuous,
\item $f$ has {$L^4_{loc}$ weak horizontal derivatives},
\item $f$ is weakly contact,
\item $f$ satisfies the distortion estimate, that is, there exists a positive and finite constant $K$ such that
\begin{equation}\label{eq:dist_quot}
 \frac{\|D_H f(p)\|}{\ell[D_H f(p)]}\leq K,\quad\text{for almost every }p\in \Heis.
\end{equation}
The quotient on the left-hand side of \eqref{eq:dist_quot} is by convention set equal to $1$ if $\|D_H f(p)\|=\ell[D_H f(p)]=0$.
\end{itemize}
If \eqref{eq:dist_quot} holds for some $K$ we also say that $f$ is {\it $K$-quasiregular}.
\end{definition}

\begin{example} Let us spell out explicitly Definition
\ref{d:QR_def} for the case of the standard sub-Riemannian Heisenberg group,
$M=\Heis$ and $g_M= g_{\Heis}$.

Assume that $f:\Heis\to \Heis$ is a  nonconstant quasiregular
according to Definition \ref{d:QR_def}. By setting $h$ (in the
definition of $L^4_{loc}$ weak horizontal derivatives) equal to a
projection on one of the coordinates, one sees that the components
of $f$ lie in $HW_{loc}^{1,4}(\Heis)$. Moreover, in this case,
$\alpha_{\Heis}(Xf)=\alpha_{\Heis}(Yf)=0$ and
\begin{displaymath}
D_H f= \begin{pmatrix} Xf_1 & Yf_1\\ Xf_2 & Yf_2\end{pmatrix}
\end{displaymath}
with respect to the basis $\{X,Y\}$.
It is well known that \eqref{eq:dist_quot} is equivalent to
\begin{displaymath}
 \|D_H f\|^4 \leq K' (\det D_H f)^2,\quad\text{a.e.}
\end{displaymath}
 for $K'=K^2$, see also Proposition \ref{p:qr_formal_jac}.

 This discussion shows that $f$ is quasiregular according to the definition commonly used for mappings in the Heisenberg group
 \cite{Da} (called ``quasiregular in the sense of Dairbekov'' in \cite{GL}).
The two definitions are in fact equivalent in this setting. In \cite{Da}, the regularity condition requires the components of the mapping to belong to the horizontal Sobolev space $HW_{loc}^{1,4}$, but this implies that $f$ has $L^4_{loc}$ weak horizontal derivatives.
 To see this, we have to verify that $h\circ f$ belongs to $HW_{loc}^{1,4}$ for all smooth functions $h:\Heis \to \mathbb{R}$, yet this follows from Proposition \ref{p:chain} applied to $M= \Heis$ and $\Psi_i$, $i\in \{1,2,3\}$, the $i$-th coordinate function in our model.
\end{example}

\subsection{Equivalent characterizations of
quasiregularity}\label{ss:equiv_char}

In this section, we give equivalent formulations of Definition
\ref{d:QR_def}. The first one allows
to interpret quasiregularity in charts. The second one is
essentially a reformulation of the distortion condition
\eqref{eq:dist_quot} in terms of the Jacobian.


\subsubsection{Contactomorphic coordinates}

\begin{definition}
Let $\Psi=(x,y,t)$ be a system of smooth $\Heis$-valued coordinate charts on $(M,HM,g_M)$ with the property that for every point $p\in M$, there exists a neighborhood $U$ and a coordinate function $\Psi:U \to \Psi(U)\subseteq \Heis$ for which
$$\Psi_{\ast,q}(HM)=H_{\Psi(q)}\Heis,\quad q\in U.$$%
 We call such coordinates \emph{contactomorphic}.
 \end{definition}

An example of contactomorphic coordinates is provided by
\emph{Darboux's theorem}, which allows us to arrange locally
$\Psi^{\ast} \alpha_{\Heis}= \alpha_M$ for a contact form
$\alpha_M$ on $M$; see for instance \cite[Theorem 3.1]{MR2682326}.

In Proposition \ref{p:equiv} below we will show how the
quasiregularity condition can be expressed in contactomorphic
coordinate charts. In the proof, the
following auxiliary result is used.



\begin{lemma}\label{eq:D_Hf_in_coord}
Assume that $U$ is a domain in $\Heis$. Let $f:U \to M$ be continuous with weak horizontal derivatives in $L_{loc}^4$. Assume further that $f$ is weakly contact and let $\Psi:f(U) \to \Psi(f(U))\subset \Heis$ be a Darboux chart so that $\Psi \circ f$ is weakly contact. Then, for almost every $q\in U$, one has
\begin{equation}\label{eq:Df_in_Psi}
 D_H f(q) = \begin{pmatrix}X(\Psi_1 \circ f)(q) & Y(\Psi_1 \circ f)(q)\\X(\Psi_2 \circ f)(q) & Y(\Psi_2 \circ f)(q)\end{pmatrix},
\end{equation}
where the matrix on the right is computed with respect to the bases
\begin{equation}\label{eq:basis1}
\{X_q,Y_q\}
\end{equation}
for $H_q\Heis$ and
\begin{equation}\label{eq:bases2}
\{X_{\Psi(f(q))}=\partial_{\Psi_1}+ 2
\Psi_2(f(q))\partial_{\Psi_3}, Y_{\Psi(f(q))}=\partial_{\Psi_2}- 2
\Psi_1(f(q))\partial_{\Psi_3}\}
\end{equation}
for $H_{f(q)}M$.
\end{lemma}

\begin{proof}
By the definition of the formal horizontal differential, we have
\begin{displaymath}
D_H f(q) (aX_q + b Y_q)= a Xf(q) + bYf(q)
\end{displaymath}
pointwise almost everywhere or in the sense of distributions. As the definition of $Xf$ and $Yf$  in \eqref{eq:def_Xf_Yf} is independent of the choice of coordinates, we may in particular work in the coordinates given by $\Psi$. Thus,
\begin{equation}\label{eq:aXbY}
a Xf(q) + bYf(q) = a \sum_{i=1}^3 X(\Psi_i \circ f)(q)
\partial_{\Psi_i}|_{f(q)}+b \sum_{i=1}^3 Y(\Psi_i \circ f)(q)
\partial_{\Psi_i}|_{f(q)}
\end{equation}
which in turn equals
$$
\begin{pmatrix}X(\Psi_1 \circ f)(q) & Y(\Psi_1 \circ f)(q)\\X(\Psi_2 \circ f)(q) & Y(\Psi_2 \circ f)(q)\end{pmatrix} \begin{pmatrix}a\\b\end{pmatrix}
$$
when expressed with respect to the basis $\{X_q,Y_q\}$ in the
source space and the basis \eqref{eq:bases2} in the target. Here
we have used the fact that $\Psi\circ f$ is a weakly contact map
between domains in the Heisenberg group, so that the right hand
side of \eqref{eq:aXbY} can be rewritten using the contact
equations in the Heisenberg group \cite[\S B]{KoRe1985} and the local frames.
\end{proof}

\begin{prop}\label{p:equiv}
Let $\Heis$ be the standard sub-Riemannian Heisenberg group and let $M$ be a contact sub-Riemannian $3$-manifold as above.
A map $f:\Heis \to M$ is quasiregular in the sense of Definition \ref{d:QR_def} if and only if
\begin{enumerate}
\item[(i)]\label{i} at every point there exists a contactomorphic
coordinate chart $\Psi$ such that $\Psi \circ f$ is quasiregular
between domains in $\Heis$ with its standard
sub-Rie{\-}mannian structure, and
\item[(ii)]\label{ii} the coordinate-free distortion estimate \eqref{eq:dist_quot} holds for a fixed constant $K$.
\end{enumerate}
\end{prop}

If the first condition holds for one contactomorphic chart $\Psi$
at $p$, it holds in fact for all such charts $\Phi$.
Moreover, assuming condition (i) of the proposition, we have that $\Psi_i \circ f$ has
weak horizontal derivatives in $L_{loc}^4$ for $i\in \{1,2,3\}$.
Then the same holds true for $h \circ f$ where $h$ is an arbitrary
smooth function by Proposition \ref{p:chain}.
Thus $Xf$ and $Yf$ can be defined, and
condition (ii) makes sense in this situation.

\begin{proof}
First, assume that $f$ satisfies conditions (i) and (ii) in Proposition \ref{p:equiv}. Then it is continuous,
has weak horizontal derivatives in $L_{loc}^{4}$ and is weakly contact. The distortion estimate
\eqref{eq:dist_quot} holds by assumption. This proves one implication.

Second, suppose that $f$ is quasiregular in the sense of
Definition \ref{d:QR_def}. Then $f$ already satisfies
condition (ii) in Proposition \ref{p:equiv} and it suffices to check condition (i). In a
neighborhood $V$ of every point $p\in M$, we can consider a
contactomorphic chart $\Psi: V \to \Heis$ given by
Darboux's theorem. The map $\Psi \circ f: U \to \Heis$, for
$U\subseteq \Heis$ small enough so that $f(U) \subseteq V$,
is continuous. By definition, $\Psi_i \circ f\in
HW_{loc}^{1,4}(U)$ for $i\in \{1,2,3\}$. Moreover, the weak
contact condition of $f$ implies that $\Psi\circ f$ is weakly contact with respect to the standard structure in source and target.

Having established the expression of the formal horizontal differential in Darboux coordinates in Lemma \ref{eq:D_Hf_in_coord}, we proceed with the proof.
We may choose $U$ small enough such that there exist constants $c,C>0$ so that
\begin{equation}\label{eq:norm_comp}
c g_{\Heis,q'}(v,v) \leq (\Psi_{\ast} g_M)_{q'}(v,v)
\leq C g_{\Heis,q'}(v,v)
\end{equation}
for all $q'\in \Psi(f(U))$ and $v \in H_{q'} \Heis$.

By assumption, the distortion estimate
\begin{displaymath}
\frac{\|Df(q)\|}{\ell[Df(q)]}\le K
\end{displaymath}
holds for a.e.\ $q\in U$, where $\|\cdot\|$ and $\ell[\cdot]$ are computed with respect to
the metric $g_M$ in the target (which corresponds to $(\Psi_{\ast})g_M$ if
$D_H f(q)$ is expressed in coordinates as in \eqref{eq:Df_in_Psi}).
The inequalities in \eqref{eq:norm_comp} allow us to switch to the norm $g_{\Heis}$ in the target.
We conclude that $\Psi \circ f$ satisfies the distortion estimate \eqref{eq:dist_quot} in a neighborhood of $p$ for {some}
constant $K'$  (depending on $c$ and $C$).
\end{proof}

\begin{remark}\label{r:Lusin}
Several deep properties of nonconstant quasiregular mappings (such as
discreteness, openness, Lusin property and vanishing measure of
the branch set) follow trivially from Proposition \ref{p:equiv} by
expressing the mappings in charts and by relying on the rich
theory in the Heisenberg group as developed in \cite{MR1778673}
and \cite{MR1630785}. It also follows that quasiregular maps
are differentiable almost everywhere in the sense of \cite{GNW,GL}
with the formal horizontal differential almost everywhere equal to
the restriction of the Margulis-Mostow derivative; cf.\ the
discussion in \cite[\S 3.2]{CLDO}.
\end{remark}

\subsubsection{Distortion estimate in terms of the Jacobian}\label{ss:JacQR}

The distortion estimate \eqref{eq:dist_quot} in Definition
\ref{d:QR_def} can be reformulated in terms of a {Jacobian
determinant} of $f$.

\begin{definition} For almost every $p\in \Heis$, the \emph{formal (horizontal) Jacobian determinant} $\det D_H f(p)$ of a map $f:\Heis \to M$ with weak horizontal derivatives is defined  as the determinant of the matrix representation of $D_H f(p)$ with respect to
 the orthonormal bases $\{X_p,Y_p\}$ and $\{e_1,e_2\}$ for $(H_p \Heis,g_{\Heis})$ and $(H_{f(p)} M, g_M)$.
\end{definition}

Note that the sign of $\det D_H f(p)$ depends on the orientation of $\{e_1,e_2\}$, but this is irrelevant for  the following proposition.

\begin{prop}\label{p:LinAlg}
Let $f:\Heis \to (M,HM,g_M)$ be a weakly contact map and, for $p\in \Heis$, let $\{e_1,e_2\}$  be a local orthonormal frame on $HM$ around $f(p)\in M$. Then, with respect to the bases $\{X_p,Y_p\}$ and $\{e_{1,f(p)},e_{2,f(p)}\}$,
\begin{equation}\label{eq:D_Hf}
D_H f(p) = \begin{pmatrix} g_M(Xf,e_1) & g_M(Yf,e_1) \\ g_M(Xf,e_2) & g_M(Yf,e_2),\end{pmatrix}
\end{equation}
and
\begin{equation}\label{eq:det}
\|D_H f(p)\| \ell[D_H f(p)] = |\det D_H f(p)|.
\end{equation}
\end{prop}

\begin{proof}
The expression \eqref{eq:D_Hf} is immediate if one expands the
vectors $Xf$ and $Yf$ in $\{e_1,e_2\}$ as
\begin{displaymath}
 Xf = g_M(Xf,e_1) e_1 + g_M(Xf,e_2) e_2\quad\text{and}\quad Yf = g_M(Yf,e_1) e_1 + g_M(Yf,e_2) e_2.
\end{displaymath}
Once $D_H f$ is expressed as a matrix with respect to the basis
$\{X,Y\}$ in the source and the basis $\{e_1,e_2\}$ in the target,
the identity \eqref{eq:det} becomes a standard fact from linear
algebra; one simply has to observe that the eigenvalues of the
symmetric matrix $(D_H f)^T (D_H f)$ are $\|D_H f\|^2$ and
$\ell[D_H f]^2$.
\end{proof}

Proposition \ref{p:LinAlg} yields the following characterization:

\begin{prop}\label{p:qr_formal_jac}
Let $M$ be a smooth $3$-manifold endowed with an equiregular distribution $HM$ and a sub-Riemannian metric $g_M$. A map $f:\Heis \to M$ is $K$-quasiregular if and only if
\begin{itemize}
\item $f$ is continuous,
\item $f$ has {$L^4_{loc}$ weak horizontal derivatives},
\item $f$ is weakly contact,
\item  the estimate
\begin{equation}\label{eq:dist}
\|D_H f\|^4 \leq K^2 (\det D_H f)^2
\end{equation}
holds almost everywhere.
\end{itemize}
\end{prop}

The formal horizontal Jacobian is related to the usual Jacobian (with respect to Popp measure) if the map is smooth. For a diffeomorphism $\phi$ between domains in $\Heis$ and a contact $3$-manifold $N$, endowed with the contact forms $\alpha_{\Heis}$ and $\alpha_N$ respectively, the Jacobian $J_{\phi}$ is given by the equation
\begin{displaymath}
\phi^{\ast} (\alpha_N\wedge \mathrm{d}\alpha_N) = J_{\phi} \,
\alpha_{\Heis} \wedge \mathrm{d}\alpha_{\Heis}.
\end{displaymath}

\begin{prop}\label{p:Jac2}
Let $V$ be a domain in $\Heis$ and $V'$ a domain in $N$. Assume that $\phi: V \to V'$ is a smooth contact transformation. Then $J_{\phi} =(\det D_H \phi)^2$. 
\end{prop}

\begin{proof}
Since $\phi$ is smooth we may apply the usual calculus for differential forms to obtain
\begin{align*}
J_{\phi} \, \alpha_{\Heis}\wedge \mathrm{d}\alpha_{\Heis}&=
\phi^{\ast} (\alpha_N \wedge \mathrm{d}\alpha_N)= \lambda^2 \,
\alpha_{\Heis}\wedge \mathrm{d}\alpha_{\Heis},
\end{align*}
where $\phi^{\ast}\alpha_N = \lambda \alpha_{\Heis}$.
Hence
$
J_\phi = \lambda^2.
$
On the other hand, if $\{e_1,e_2\}$ is a local orthonormal frame on $V'$, and $\{\nu_1,\nu_2\}$ is the dual frame of $1$-forms, then $X\phi= \nu_1(X\phi) e_1 + \nu_2(X\phi) e_2$ and $Y\phi = \nu_1 (Y\phi) e_1 + \nu_2(Y\phi) e_2$, whence
\begin{displaymath}
D_H \phi = \begin{pmatrix} \nu_1(X\phi)& \nu_1(Y\phi)\\ \nu_2(X\phi)&\nu_2(Y\phi) \end{pmatrix}
\end{displaymath}
with respect to the bases $\{X,Y\}$ and $\{e_1,e_2\}$. Since
$
X\phi = \phi_{\ast} X$ and  $Y\phi= \phi_{\ast}Y$,
it follows that
\begin{equation*}\begin{split}
\det D_H \phi
&= \nu_1(\phi_{\ast}X)\nu_2(\phi_{\ast}Y) -\nu_1(\phi_{\ast}Y)\nu_2(\phi_{\ast}X) \\
&=\mathrm{d}\alpha_N(\phi_{\ast}X,\phi_{\ast}Y)
= \lambda \mathrm{d}\alpha_{\Heis}(X,Y)
=\lambda.
\end{split}\end{equation*}
The claim follows.
\end{proof}

\section{Proof of the main theorem}\label{s:main}
We now prove Theorem \ref{thm:main}, following the steps outlined in the introduction. We fix an equiregular sub-Riemannian $3$-manifold $M$ whose fundamental group has growth rate larger than $4$.

\subsection{Topology and covering theory}
\label{sec:topology}

By definition, there exists a finitely generated subgroup $\Gamma$ of $\pi_1(M)$ so that $\Gamma$ has growth rate larger than $d$, for some number $d>4$. We will associate to $\Gamma$ a relatively compact ``core'' $M'\subset M$. The goal of this section is to prove that a specific lift $\widetilde{M}''$ of $M'$ inside the universal cover $\widetilde{M}$ of $M$ is quasi-isometric to $\Gamma$.

Fix a basepoint $x_0$ in $M$ and smooth closed curves $\gamma_1, \ldots, \gamma_s$ for some $s\geq 2$ such that $\Gamma$ is generated by $[\gamma_1], \ldots, [\gamma_s]\in \pi_1(M)$. We may assume, without loss of generality, that the curves are simple, intersect only at the basepoint, and intersect transversally at the basepoint. Let $\overline{M'} $ be a closed, connected manifold with $\mathcal{C}^{\infty}$ boundary such that $M':=\mathrm{int}(\overline{M'})$ satisfies:
\begin{enumerate}
\item $\mathrm{int}(M') = M' \subset \overline{M'} \Subset M$,
\item  $\gamma_i \subset M'$ for all $i\in \{1,\ldots,s\}$,
\item $\pi_1(M')$ is the free group generated by $[\gamma_1],\ldots,[\gamma_s]$.
\end{enumerate}

Let $\widetilde M$ and $\widetilde{M'}$ denote the universal covers of $M$ and $M'$ respectively. While $\widetilde {M'}$ is (generically) not a subset of $\widetilde M$, there is some intermediate cover of $M'$ inside $\widetilde M$. Indeed, the inclusion of $M'$ in $M$ induces a map $\widetilde{M'} \to \widetilde{M}$ whose image we denote by $\widetilde{M}''$.

\begin{lemma}
\label{l:PropGroupAction}
The following properties hold:
\begin{enumerate}
\item $\widetilde{M}''$ is a cover of $M'$ under the standard
projection $\pi:\widetilde{M} \to M$, \item the action of $\Gamma$
on $\widetilde{M}$ leaves $\widetilde{M}''$ invariant, \item
$\widetilde{M}''/\Gamma = M'$.
\end{enumerate}
All these properties also hold if $M'$ is replaced by its closure
$\overline{M'}$. That is, one
can define in the same way a cover space of $\overline{M'}$ so
that $\Gamma$ acts on this cover, and the quotient of the action
can be identified with $\overline{M'}$.
\end{lemma}

The universal cover $\widetilde{M}$ can be endowed with a contact sub-Riemannian structure by lifting the contact form and metric from $M$. Let us denote by $d_M$, resp.\ $d_{\widetilde{M}}$, the sub-Riemannian metric on $M$, resp.\ $\widetilde{M}$.
The embeddings $M'\hookrightarrow M$ and $\widetilde{M}'' \hookrightarrow \widetilde{M}$ equip $M'$ and $\widetilde{M}''$ with sub-Riemannian metrics (denoted $\delta_{M'}$ and $\delta_{\widetilde{M}''}$) such that the covering map $\pi:\widetilde{M}'' \to M'$ becomes a local isometry. For example, the distance between two points of $M'$ is the infimal $g_M$-length of horizontal curves contained in $M'$ joining these two points. We call these quantities the {\it intrinsic distance} on $M'$ and $\widetilde{M}''$.

Understanding the topological and metric properties of a
submanifold endowed with an \emph{intrinsic} distance
is more
challenging in the present sub-Riemannian setting than it would be
in the Riemannian case. The difficulties arise already in case $M$
is the Heisenberg group $\Heis$ itself. For instance, Monti and
Rickly showed in \cite{MR2135806} that the only geodetically
convex subsets of $\Heis$ are the empty set, points,
geodesic arcs, and the whole space. Moreover, there exist domains
in $\Heis$, even $\mathcal{C}^1$-smooth ones, for which some
points on the boundary cannot be joined from inside the domain by
rectifiable curves \cite{MR2021252}. These complications hint at
the subtleties involved in analyzing the intrinsic distance on
submanifolds with boundary in a sub-Riemannian manifold. In the
present section, we discuss properties of the intrinsic distances
on $M'$ and $\widetilde{M}''$.

First note that although $\overline{M'}$ is compact in the topology of $M$, the intrinsic distance $\delta_{M'}$ might induce a different topology on $M'$ and we do not a priori know whether $M'$ is bounded with respect to this distance. Our first goal is to show that $\delta_{M'}$ is bi-Lipschitz equivalent to the \emph{extrinsic} distance $d_M|_{M'}$, in a way that extends to the boundary of $M'$.

We define the function
\begin{equation}\label{d-intr}
d_{\mathrm{intr}}:\overline{M'}\times \overline{M'} \to [0,+\infty],\quad d_{\mathrm{intr}}(p,q)=\inf \mathrm{length}(\gamma),
\end{equation}
where the infimum is taken over curves $\gamma:[0,1] \to \overline{M'}$ such that $\gamma(0)=p$, $\gamma(1)=q$, and most importantly, $\gamma(x)\in M'$ for $x\in (0,1)$. The length in \eqref{d-intr} is computed with respect to the sub-Riemannian metric tensor $g_M$, so it agrees with the length in the extrinsic metric $d_M|_{M'}$. In particular, non-horizontal curves have infinite length. It is clear that the restriction of $d_{\mathrm{intr}}$ to $M'$ agrees with $\delta_{M'}$.

A priori, the value of $d_{\mathrm{intr}}(p,q)$ could be infinite. The following proposition shows that this is not the case.

\begin{prop}\label{p:biLip_metric}
The function $d_{\mathrm{intr}}$ defines a metric on $\overline{M'}$ that is bi-Lipschitz equivalent to $d_M|_{\overline{M'}}$.
\end{prop}

\begin{proof}
It is immediate that
\begin{displaymath}
d_{\mathrm{intr}}(p,q)\geq d_M(p,q),\quad \text{for all }p,q\in \overline{M'}.
\end{displaymath}
Since $d_M$ is a metric, this shows that $d_{\mathrm{intr}}$ is non-degenerate, and it suffices to prove that also the reverse inequality holds -- up to a multiplicative constant.

We first show that the two distances are locally bi-Lipschitz equivalent. This implies that $d_{\mathrm{intr}}$ induces the original topology on $\overline{M'}$. We conclude that $(\overline{M'},d_{\mathrm{intr}})$ is compact and use this information to show that the two metrics are globally bi-Lipschitz equivalent. A similar argument can be found in \cite[\S3]{MR3226622}.

We fix a point $p\in \overline{M'}$ and a Darboux chart $\Phi: U \to V$ mapping a neighborhood $U$ of $0$ in $\Heis$ to a neighborhood $V$ of $p=\Phi(0)$ in $M$.
This is possible since $M$ is a contact manifold. The chart map $\Phi$ pulls the metric tensor $g_M$ back to a sub-Riemannian metric $g_U$ on $U$, which, on compact sets, is comparable to the standard sub-Riemannian metric $g_{\Heis}$.

Furthermore, there exist relatively compact neighborhoods $U_0$ of
$0$ in $U$, and $V_0$ of $p$ in $V$, with $\Phi(U_0)=V_0$, such
that the $d_{\Heis}$-distance between points of $U_0$ is
realized by curves contained in $U$, while the $d_M$-distance
between points in $V_0$ is realized by curves contained in $V$.
Thus, in order to prove that $d_{\mathrm{intr}}$ and
$d_M|_{\overline{M'}}$ are bi-Lipschitz equivalent on
$V_0\cap \overline{M'}$, it suffices to prove that there exists a
finite constant $c_p$ such that
\begin{equation}\label{eq:goal_BL_metric}
\inf_{\substack{\gamma \\ \gamma((0,1)) \subset \Phi^{-1}(V\cap
M')}} \mathrm{length}_{g_{\Heis}}(\gamma) \leq c_p
\inf_{\gamma} \mathrm{length}_{g_{\Heis}}(\gamma)
\end{equation}
for all $u,u'\in  \Phi^{-1}(V_0 \cap \overline{M'})$ and for
$\gamma:[0,1] \to U$, with $\gamma(0)=u$ and $\gamma(1)=u'$.

First, suppose that $p \in M'$. In this case, by choosing $U$
smaller if necessary, we may assume that $\Phi(U) \cap \partial M'
= \emptyset$, whence $\Phi^{-1}( V\cap M')= U$ and
\eqref{eq:goal_BL_metric} clearly holds.

Now suppose that $p \in \partial \overline{M'}$. By making $U$
smaller if necessary, we may assume that $\Phi^{-1}(V\cap \partial
M')$ is a smoothly embedded disk which separates $U$ in two
domains.  Denote by $U'$ the domain $\Phi^{-1}(V\cap M')$ in
$U$. We may further assume that the subdomain $U_0$ is chosen so
that $\Omega := \Phi^{-1}(V_0 \cap M') \subset U_0$ is a domain
with $\mathcal{C}^{1,1}$ boundary with $0 \in \partial\Omega$ and
\begin{displaymath}
\partial \Omega = \left(\Phi^{-1}(V_0 \cap \partial M') \right) \dot{\cup} \left(\partial \Omega \cap U'\right).
\end{displaymath}
By \cite[Theorem 1.3]{MR2135732}, $\Omega$ is an NTA domain in $(\Heis,d_{\Heis})$ and hence also a uniform domain (see, for instance, \cite[Proposition 4.2]{MR1323792}). By a limiting argument (see the remark on p.\ 270 of \cite{MR1323792}), it follows that points in $\overline{\Omega}$ can be joined by uniform curves. In particular, $\overline{\Omega}$ is quasiconvex: there exists $C>0$ so that for any pair of points $u,u'\in \overline\Omega$ there exists a rectifiable curve $\gamma:[0,1] \to \overline{\Omega}$ such that $\gamma((0,1))\subset \Omega$, $\gamma(0)=u$, $\gamma(1)=u'$, and $\mathrm{length}_{g_{\Heis}}(\gamma) \leq C d_{\Heis}(u,u')$. Since $u$ and $u'$ lie in $U_0$, to compute $d_{\Heis}(u,u')$ it suffices to consider curves contained in $U$. Thus we have established \eqref{eq:goal_BL_metric} when $p \in \partial M'$. We conclude that for every point $p\in \overline{M'}$ there exists an open neighborhood $V_p$ of $p$ in $M$ and a finite constant $C_p>0$ such that
\begin{displaymath}
d_M(q,q') \leq d_{\mathrm{intr}}(q,q') \leq C_p d_M(q,q'),\quad \text{for all }q,q'\in V_p \cap \overline{M'}.
\end{displaymath}
In particular, the intrinsic and extrinsic topologies on $\overline{M'}$ agree. Hence $\overline{M'}$ is compact and therefore bounded with respect to the metric $d_{\mathrm{intr}}$.

Since $\overline{M'}$ is compact, we can cover $\overline{M'}$ by finitely many open sets $V_{p_1},\ldots, V_{p_N}$ of the above form. By the Lebesgue number lemma, there exists $r_0>0$ such that for all $0<r<r_0$ and $p\in \overline{M'}$, we have
\begin{equation}\label{eq:LebesgueNumber}
\{q\in \overline{M'}: d_M(p,q) < r\}\subseteq V_{p_i}\quad \text{for some }i\in \{1,\ldots,N \}.
\end{equation}
In order to prove that $d_{\mathrm{intr}}$ and $d_M$ are bi-Lipschitz equivalent on $\overline{M'}$, we have to show that
\begin{displaymath}
\sup_{p\neq q}\frac{d_{\mathrm{intr}}(p,q)}{d_M(p,q)}
\end{displaymath}
is uniformly bounded. Considering the cases $d_M(p,q)<r$ and $d_M(p,q)\ge r$ in turn, we see that
$$
\frac{d_{\mathrm{intr}}(p,q)}{d_M(p,q)} \le \max \left\{ C_{p_1},\ldots,C_{p_N},\frac{\mathrm{diam}_{d_{\mathrm{intr}}}(\overline{M'})}{r} \right\}.
$$
for all $p \ne q$. This completes the proof.
\end{proof}

The restriction of the metric $d_{\mathrm{intr}}$ to $M'$ is by definition a length metric on $M'$ and agrees with the sub-Riemannian distance $\delta_{M'}$ induced by $g_M$. Using the covering map $\pi$, this distance can be lifted to a length metric $d_{\widetilde{M}''}$ on $\widetilde{M}''$ which agrees with the sub-Riemannian distance $\delta_{\widetilde{M}''}$ induced on $\widetilde{M}''$ by the pull-back of $g_M$ under $\pi$.

\begin{lemma}\label{l:GammaQI}
The metric space $(\widetilde{M}'',d_{\widetilde{M}''})$ is quasi-isometric to $\Gamma$ endowed with a word metric.
\end{lemma}

Recall that a map $f:X\to Y$ between metric spaces is {\it $(A,B)$-quasi-isometric} for $A \ge 1$ and $B>0$ if $A^{-1}d_X(x_1,x_2)-B \le d_Y(f(x_1),f(x_2)) \le A d_X(x_1,x_2)+B$ for all $x_1,x_2 \in X$ and if $f(X)$ is $B$-coarsely dense in $Y$, i.e., every point of $Y$ is within distance $B$ of $f(X)$.

\begin{proof}
The set $\overline{M'}$ can be endowed with a length metric $d_{\overline{M'}}$, where the distance between two points $p,q$ is defined by minimizing the $g_M$-length of curves in $\overline{M'}$ connecting $p$ and $q$. This distance is bounded from below by the restriction of $d_M$ to $M'$, and bounded from above by $d_{\mathrm{intr}}$. By Proposition
\ref{p:biLip_metric}, $d_{\overline{M'}}$ is comparable to both $d_M|_{\overline{M'}}$ and $d_{\mathrm{intr}}$. The metric $d_{\overline{M'}}$ lifts to a length metric $d_N$ on the component $N$ of $\pi^{-1}(\overline{M'})$ that contains $\widetilde{M}''$.

Since $(N,d_N)$ and $(\widetilde{M}'',d_{\widetilde{M}''})$ are quasi-isometric, it suffices to prove that $(N,d_N)$ is quasi-isometric to $\Gamma$.
This follows from the Milnor--\v{S}varc lemma, upon observing that $(N,d_N)$ is a length space on which $\Gamma$ acts properly discontinuously and cocompactly by isometries. Here we have used Lemma \ref{l:PropGroupAction} and the fact that $\overline{M'}$ is compact in the metric $d_{\overline{M'}}$, which follows from Proposition \ref{p:biLip_metric}.
\end{proof}

\subsection{Proof of the rough isoperimetric inequality}\label{ss:isop}

We will show in Section \ref{sec:globalisoperimetric} that $\widetilde{M}$ satisfies a $d$-dimensional isoperimetric inequality for some $d>4$. To this end, we first establish in Proposition \ref{p:roughIP} a rough, or coarse,  $d$-dimensional isoperimetric inequality  for a net on $\widetilde{M}''$.

Following the terminology used by Kanai in \cite{MR792983}, we call a countable set $Y$ a \emph{net} if there exists a set-valued function
\begin{displaymath}
N: Y \to \{A:\, A \subseteq Y\},\quad y\mapsto N(y)
\end{displaymath}
with the property that
\begin{enumerate}
\item $N(y) \subseteq Y$ is finite,
\item $x\in N(y)$ if and only if $y\in N(x)$.
\end{enumerate}
The points in $N(y)$ are called the \emph{neighbors} of $y$, and $Y$ is said to be \emph{connected} if for any two points $y$ and $x$ in $Y$ there exists a chain of finitely many points connecting $y$ and $x$ so that any two consecutive points are neighbors. The \emph{combinatorial distance} $\delta(y,x)$ is the minimal length (number of elements in the chain) a path must have to connect $y$ and $x$.

The nets considered in this paper are all connected, and one might think of them as {graphs}. We encounter two types of nets:
\begin{enumerate}
\item finitely generated groups, where the combinatorial metric agrees with the word metric with respect to the system of generators,
\item point sets $Y$ on a metric space $(X,d)$ that are \emph{$\varepsilon$-separated} ($d(y,x)\geq \varepsilon$ for all $y,x\in Y$) and \emph{maximal} (with respect to order of inclusion), and where $N(y)=\{x\in P: 0<d(y,x)\leq 2\varepsilon\}$.
\end{enumerate}

We recall that in every metric space there exist maximal $\varepsilon$-separated nets for every $\varepsilon>0$. Moreover, a totally bounded metric space contains a \emph{finite} maximal $\varepsilon$-separated net for every $\varepsilon>0$. In particular, every compact set in a metric space contains a finite maximal $\varepsilon$-separated net for every $\varepsilon>0$.

A net $N$ is \emph{uniform} if
\begin{displaymath}
\sup_{y\in Y} \sharp N(y) <\infty.
\end{displaymath}
Nets derived from finitely generated groups in the above way are always uniform. The metric space $(\widetilde{M}'',d_{\widetilde{M}''})$ defined at the beginning of this section admits a uniform net, which is moreover maximally $\varepsilon$-separated for some $\varepsilon>0$.

\begin{lemma}\label{l:unif}
There exists $\varepsilon_0>0$ such that for every $0<\varepsilon<\varepsilon_0$, the manifold $\widetilde{M}''$ contains a maximal $\varepsilon$-separated net $Y$ such that
\begin{equation}\label{eq:UDBG}
\sup_{y\in Y} \sharp \left(B_{\widetilde{M}''}(y,r)\cap Y\right) <\infty,\quad\text{for all }0<r<\infty.
\end{equation}
In particular,  $Y$ is uniform.
\end{lemma}

\begin{proof}
Since $\pi|_{\widetilde{M}''}:\widetilde{M}''\to M'$ is a locally isometric cover, and since $\overline{M'}$ is compact, there exists a constant $\varepsilon_0>0$ such that $\pi$ is an isometry when restricted to $\varepsilon_0$-balls in $\widetilde{M}''$. We let $\varepsilon$ be a positive number less than $\varepsilon_0$ and pick a maximal $\varepsilon$-separated net $Y_{M'}$ in $M'$. Notice that this net is finite since $\overline{M'}$ is compact.

We denote the lift of $Y_{M'}$ to $\widetilde{M}''$ by $Y$, that is,
\begin{displaymath}
Y=\{\gamma . x: \gamma \in \Gamma,\, x\in (\pi|_{\widetilde{M}''})^{-1}(y),\,y\in Y_{M'}\}.
\end{displaymath}
Clearly, $Y$ is $\varepsilon$-separated. It is also maximal, for if we could add a point $p\in \widetilde{M}''\setminus Y$ with $d_{\widetilde{M}''}(p,y)\geq \varepsilon$ for all $y\in Y$, then $\pi(p)$ would be at distance at least $\varepsilon$ from every point in $Y_{M'}$, which contradicts the maximality of $Y_{M'}$. We observe further that, by construction, $Y$ is $\Gamma$-invariant.

Denote by $n$ the cardinality of $Y_{M'}$. We choose points $y_1,\ldots,y_n$ in $Y$ such that
\begin{displaymath}
\pi \{y_1,\ldots,y_n\}=Y_{M'}.
\end{displaymath}
Each $y\in Y$ can be written as $y =\gamma . y_i$ for a unique $i\in \{1,\ldots,n\}$ and a unique $\gamma \in \Gamma$. Thus
\begin{displaymath}
d_{\widetilde{M}''}(y,\gamma .y_1)\leq \max_{1\leq i\leq n} d_{\widetilde{M}''}(y_i,y_1).
\end{displaymath}
This shows that  there exists a $n$-to-$1$ quasi-isometry from $(Y,d_{\widetilde{M}''})$ to
\begin{displaymath}
Y_0 = \{\gamma. y_1: \gamma \in \Gamma\} \quad\text{with}\quad d_{\widetilde{M}''}
\end{displaymath}
given by $y = \gamma. y_i \mapsto \gamma. y_1$. On the other hand,
we know from the proof of the Milnor--\v{S}varc lemma, see for
instance \cite[Theorem 1.18]{MR2007488}, that
$(Y_0,d_{\widetilde{M}''})$ is quasi-isometric to $\Gamma$ with
the word metric via the map $\gamma.y_1 \mapsto \gamma$. It
follows that there is an $n$-to-$1$ quasi-isometry $\varphi: Y \to
\Gamma$. Hence, for every $r>0$, there exists $r'$, depending on
$r$ and the quasi-isometry constants of $\varphi$, such that
\begin{displaymath}
\varphi \left(B_{\widetilde{M}''}(y,r)\cap Y \right) \subset B_{\Gamma}(\varphi(y),r'),\quad \text{for all }y\in Y.
\end{displaymath}
The ball on the right hand side contains only a finite set of elements in $\Gamma$, whose cardinality can be bounded  depending on $r'$, but independently of $y$. It follows that $B_{\widetilde{M}''}(y,r)\cap Y$ contains at most this number times $n$ net points. This proves \eqref{eq:UDBG}.
\end{proof}

\begin{remark}\label{rem:net_point_overlap}
Condition \eqref{eq:UDBG} in Lemma \ref{l:unif} can be extended to all points in $\widetilde{M}''$. Indeed, if $x$ is a point in $\widetilde{M}''\setminus Y$, then by maximality of $Y$, there exists $y\in Y$ with $d_{\widetilde{M}''}(x,y)<\varepsilon$. Thus, for all $r>0$, we have that
\begin{displaymath}
B_{\widetilde{M}''}(x,r) \subset B_{\widetilde{M}''}(y,r+\varepsilon).
\end{displaymath}
It follows by \eqref{eq:UDBG} that for all $x\in \widetilde{M}''$, we have
\begin{equation}\label{eq:Kanai(2.4)}
\sharp \left(B_{\widetilde{M}''}(x,r) \cap Y \right) \leq \sup_{y\in Y} \sharp \left(B_{\widetilde{M}''}(y,r+\varepsilon)\cap Y\right)=:\nu(r,\varepsilon)<\infty.
\end{equation}
\end{remark}

The net $Y$ from Lemma \ref{l:unif} is quasi-isometric to $(\widetilde{M}'',d_{\widetilde{M}''})$ if it is seen as a subset of $\widetilde{M}''$ and endowed with $d_{\widetilde{M}''}$. The same holds true, but is less immediate, if $Y$ is equipped with the \emph{combinatorial distance} $\delta$. Clearly
\begin{displaymath}
d_{\widetilde{M}''}(y,y') \leq 2 \varepsilon \delta(y,y')
\end{displaymath}
for all $y,y'\in Y$.
One can use property \eqref{eq:Kanai(2.4)} to prove that $\delta$ is  controlled also from above in terms of $d_{\widetilde{M}''}$ (up to multiplicative and additive constants). An analogous statement is known for  complete Riemannian manifolds with lower Ricci curvature bound, and an inspection of \cite[Lemma 2.5]{MR792983} shows that it carries over to length spaces satisfying condition \eqref{eq:Kanai(2.4)}.

Let $Y$ be a net endowed with the combinatorial distance $\delta$. We define the \emph{boundary} of a set $S\subseteq Y$ as
\begin{displaymath}
\partial S = \{y\in Y:\; \delta(y,S)=1 \}.
\end{displaymath}

\begin{definition}
We say that $Y$ satisfies a \emph{rough $d$-dimensional isoperimetric inequality} if there exists a constant $0<C<\infty$ such that
\begin{displaymath}
(\sharp S)^{\frac{d-1}{d}}\leq C \,\sharp \partial S
\end{displaymath}
for all nonempty finite subsets $S$ of $Y$.
\end{definition}

\begin{prop}\label{p:roughIP}
There exists $\varepsilon_0>0$ such that for every $0<\varepsilon<\varepsilon_0$, the manifold $\widetilde{M}''$ contains a maximal $\varepsilon$-separated net $Y$ which satisfies a rough $d$-dimensional isoperimetric inequality for some $d>4$.
\end{prop}

\begin{proof}
The first part of the proof is not specific to the sub-Riemannian setting, but instead proceeds exactly the same way as \cite[Proof of Theorem 1.3]{MR2832708}. Namely one uses the fact that $\Gamma$ has a growth rate $d>4$ in order to deduce by \cite[Th\'{e}or\`{e}me 1]{MR1232845} that it satisfies a rough $d$-dimensional isoperimetric inequality.

Then we choose $\varepsilon_0$ and $Y$ as in Lemma \ref{l:unif}. By the comment after \eqref{eq:Kanai(2.4)}, we know that $(Y,\delta)$ is quasi-isometric to $(\widetilde{M}'',d_{\widetilde{M}''})$, and thus, by Lemma \ref{l:GammaQI}, also quasi-isometric to $\Gamma$ endowed with the word metric.

Kanai has shown in \cite[Lemma 4.2]{MR792983} that the validity of a rough $d$-dimensional isoperimetric inequality is a quasi-invariant for uniform nets. Recall that $\Gamma$ is uniform (since the group is finitely generated), and that $Y$ is uniform by Lemma \ref{l:unif}. It follows that $Y$ satisfies a rough $d$-dimensional isoperimetric inequality.
\end{proof}

We will later apply Proposition \ref{p:roughIP} to prove a (smooth) isoperimetric inequality on $\widetilde{M}''$.
A close inspection of Kanai's proof in the Riemannian setting reveals that the full strength of volume comparison geometry is not needed. For our purposes, a much weaker estimate suffices.

\begin{lemma}\label{eq:volume_comparison}
There exists $\varepsilon_1>0$ such that for every $\varepsilon$-net $Y$ as above with $\varepsilon<\varepsilon_1$,
there are constants $0<c_- \leq c_+ <\infty$ such that
\begin{displaymath}
c_- \leq \mu_{\widetilde{M}}(B_{\widetilde{M}''}(y,\varepsilon))\leq  c_+,\quad \text{for all }y\in Y.
\end{displaymath}
The constants $c_-$ and $c_+$ may depend on the data of the manifold and on $\varepsilon$, but not on $y$.
\end{lemma}

\begin{proof}
Recall that $Y= \{\gamma .y_i:\gamma \in \Gamma, i=1,\ldots,n\}$ for a finite set of points $\{y_1,\ldots,y_n\} \in \widetilde{M}''$. Since $\gamma$ acts by isometries and $\mu_{\widetilde{M}}$ agrees up to a multiplicative constant with the $4$-dimensional Hausdorff measure with respect to $d_{\widetilde{M}''}$, it suffices to consider the mass of the balls $B_{\widetilde{M}''}(y_i,\varepsilon)$ for $i=1,\ldots, n$, more precisely, to prove that this volume is positive and finite. If we choose $\varepsilon_1$ no larger than the constant $r_0$ from the proof of Proposition \ref{p:biLip_metric} and the constant $\varepsilon_0$ from the proof of Lemma \ref{l:unif}, then every such ball is isometric to a ball in $M'$ which is contained in one of finitely many sets that cover $M'$ and that can be mapped bi-Lipschitzly onto a domain in the Heisenberg group by a Darboux chart. It is well known that $\mu_{\Heis}(B_{\Heis}(p,r))=c r^4$ for all $r>0$, $p\in\Heis$, and a positive and finite constant $c$. The claim follows.
\end{proof}

This concludes our discussion of the coarse geometry of $\widetilde{M}''$. In the next section we will focus on the local geometry.

\subsection{Proof of the relative isoperimetric inequality}
\label{sec:localisoperimetric}

In order to derive from a rough isoperimetric inequality a smooth, global one, we need local information on the geometry of the manifold  $\widetilde{M}''$. In this section, we prove a weak Sobolev--Poincar\'{e} inequality and a weak relative $4$-dimensional isoperimetric inequality for balls of a fixed radius centered in the net $Y$. In the continuous version of the isoperimetric inequality, the cardinality of a finite set $S$ and of its boundary are replaced by the volume of a domain $\Omega$ and the perimeter of its boundary.

We now explain the notion of perimeter in our setting. We follow the presentation in \cite{MR3412382}, which is specific to sub-Riemannian contact $3$-manifolds. In \cite{MR3412382}, the manifolds are assumed to be endowed with a \emph{global} contact form, but in light of Remark \ref{r:popp_global_orientable}, this assumption is not necessary for our application.

\begin{definition}\label{eq:smoothDiv}
The \emph{divergence} of a smooth horizontal vector field $V$ on a manifold $(N,HN,g_N)$ is the real-valued function $\mathrm{div}_N V$ characterized by the identity
\begin{equation}\label{eq:Lie_deriv_div}
\mathcal{L}_V \mathrm{vol}_N = (\mathrm{div}_N V)\mathrm{vol}_N,
\end{equation}
where $\mathcal{L}_V$ denotes the Lie derivative along $V$.
\end{definition}

\begin{definition} Let $N$ be a contact sub-Riemannian $3$-manifold.  The relative
\emph{perimeter} of a measurable set $E\subset N$ in an open set $\Omega \subseteq N$ is defined to be
\begin{displaymath}
\mathcal{P}(E,\Omega) := \sup \left\{ \int_{E \cap \Omega} \mathrm{div}_N V\;\mathrm{d}\mu_N:\; V\text{ horizontal,}\, \|V\|_{\infty}\leq 1 \right\},
\end{displaymath}
where the supremum is taken over $\mathcal{C}^1$ vector fields $V$ with compact support in $\Omega$.
\end{definition}

The perimeter $\mathcal{P}(E,\Omega)$ can be seen as a measure for the area of the boundary of $E$ in $\Omega$, see for instance \cite[2.3]{MR3412382} and \cite[(2.7)]{MR3474402}.

Our first goal is to prove a local relative isoperimetric inequality on $M'$.
We follow the argument of Galli and Ritor\'{e} in \cite{MR2979606}, where such an inequality is proved for balls centred in a compact subset of a contact sub-Riemannian manifold. For our argument we only need a weak form of this relative isoperimetric inequality, namely a statement about $\varepsilon$-balls centred in the points of a given $\varepsilon$-net in $M'$. The reason why we cannot directly apply \cite[Lemma 3.7]{MR2979606}, is that this statement holds only for small enough balls, where the smallness condition depends on the compact set. In other words, if we consider balls centered in the $\varepsilon$-net $Y_{M'}$, we might only get an estimate for balls at a scale much smaller than $\varepsilon$. Since $M'$ is an open subset of $M$, this issue cannot be fixed by a simple compactness argument. We therefore give a direct proof for the result  we are going to apply. This includes analyzing intrinsic $\varepsilon$-balls in $M'$ whose closure with respect to $d_M$ might intersect the boundary of $M'$.

\begin{prop}[Weak local relative isoperimetric inequality]\label{p:localIP}
There exists $\delta_1>0$ such that for every $\varepsilon$-net $Y$, $\varepsilon< \delta_1$, on $\widetilde{M}''$, there exist constants $C_I>0$ and $1<c<\infty$, depending only on $Y$ and the data of the manifold, such that for any bounded set $E\subset \widetilde{M}''$ with finite perimeter, one has
\begin{displaymath}
C_I \left(\min \{\mu_{\widetilde{M}}(E\cap B_{\widetilde{M}''}(y,\varepsilon)),\mu_{\widetilde{M}}((\widetilde{M}'' \setminus E) \cap B_{\widetilde{M}''}(y,\varepsilon))\}\right)^{\frac{d-1}{d}}\leq \mathcal{P}(E,B_{\widetilde{M}''}(y,c\varepsilon)),
\end{displaymath}
for all $y \in Y$.
\end{prop}

The result follows from a suitable weak Sobolev-Poincar\'{e} inequality for the balls $B_{\widetilde{M}''}(y,\varepsilon)$, $y\in Y$. While we may assume that $M$ and thus also $\widetilde{M}$ are complete (by a conformal change of metric which does not affect quasiregularity), the same is not true for $M'$ and $\widetilde{M}''$
and we must prove by hand that the considered balls are John domains.

We denote by
\begin{displaymath}
u_E = \frac{1}{\mu(E)} \int_E u \;\mathrm{d}\mu= \Barint_E u \;\mathrm{d}\mu
\end{displaymath}
the mean value of a function $u:X \to \mathbb{R}$ over a measurable set $E$ with positive mass in a metric measure space $(X,d,\mu)$.
To establish the desired Poincar\'{e} inequality for $\varepsilon$-balls centred in the points of a $\varepsilon$-net $Y_M$ on $M'$, we will use Darboux charts to transfer the problem to the Heisenberg group. Balls with respect to the intrinsic distance on a domain $U$ in $\Heis$ need not be John domains even if the boundary of the ball is smooth, but they can be compared to subsets of $U$ that are John domains. In this context the following lemma is useful. The proof is a standard argument which works much more generally and which we reproduce here for completeness.

The following definition is used here and in the following:

\begin{definition}\label{d:horiz_grad}
Let $(N,HN,g_N)$ be a sub-Riemannian contact $3$-manifold. The \emph{horizontal gradient} $\nabla_H u$ of a $\mathcal{C}^1$
function $u:N \to \mathbb{R}$ is the unique horizontal vector field on $N$ with the property that
\begin{displaymath}
g_N(\nabla_H u, V) = \mathrm{d}u (V),\quad \text{for all }V\in HN.
\end{displaymath}
\end{definition}
See \cite[Section 2.2]{MR2502528} for an expression of the gradient in a local orthonormal frame.

\begin{lemma}\label{l:incl_ball_PI}
Let $\Omega$ be a domain in $\Heis$ and $B\subset \Omega$  a measurable subset with the property that
\begin{displaymath}
\left( \Barint_B |u-u_B|^{\frac{4}{3}}\;\mathrm{d}\mu_{\Heis}\right)^{\frac{3}{4}}\leq C \mathrm{diam}(B) \left(\Barint_B |\nabla_H u|\;\mathrm{d}\mu_{\mathbb{H}}\right)
\end{displaymath}
for some real-valued function $u$ defined in a neighborhood of $\Omega$.
If $B_1$ and $B_2$ are measurable sets of positive mass such that
$B_1 \subset B \subset B_2\subset \Omega,
$
then
\begin{displaymath}
\left(\Barint_{B_1} |u-u_{B_1}|^{\frac{4}{3}}\;\mathrm{d}\mu_{\Heis}\right)^{\frac{3}{4}}\leq  C' \mathrm{diam}(B_2) \left(\Barint_{B_2} |\nabla_H u|\;\mathrm{d}\mu_{\mathbb{H}}\right),
\end{displaymath}
where $C'$ depends only on $C$ and on the ratios $\mu_{\Heis}(B)/\mu_{\Heis}(B_1)$ and $\mu_{\Heis}(B_2)/\mu_{\Heis}(B)$.
\end{lemma}

\begin{proof}
For simplicity we write $\mu=\mu_{\Heis}$. First we observe that
\begin{align*}
\left(\int_{B_1} |u-u_{B_1}|^{\frac{4}{3}}\;\mathrm{d}\mu\right)^{\frac{3}{4}}& \leq
\left(\int_{B_1} |u-u_{B}|^{\frac{4}{3}}\;\mathrm{d}\mu\right)^{\frac{3}{4}}+\left(\int_{B_1} |u_{B_1}-u_{B}|^{\frac{4}{3}}\;\mathrm{d}\mu\right)^{\frac{3}{4}}\\
&= \left(\int_{B_1} |u-u_{B}|^{\frac{4}{3}}\;\mathrm{d}\mu\right)^{\frac{3}{4}}
+ \mu(B_1)^{\frac{3}{4}}\left| \Barint_{B_1} u -u_B \;\mathrm{d}\mu\right|\\
&\leq \left(\int_{B_1} |u-u_{B}|^{\frac{4}{3}}\;\mathrm{d}\mu\right)^{\frac{3}{4}}
+ \mu(B_1)^{\frac{3}{4}} \Barint_{B_1} |u -u_B| \;\mathrm{d}\mu\\
&\leq 2 \left(\int_{B_1} |u-u_B|^{\frac{4}{3}}\;\mathrm{d}\mu\right)^{\frac{3}{4}},
\end{align*}
where we have used H\"older's inequality in the last step.

With this estimate and the assumed inequality for $B$ in hand, it now follows that
\begin{align*}
\left(\Barint_{B_1} |u-u_{B_1}|^{\frac{4}{3}}\;\mathrm{d}\mu\right)^{\frac{3}{4}}&\le 2 \left(\Barint_{B_1} |u-u_B|^{\frac{4}{3}}\;\mathrm{d}\mu\right)^{\frac{3}{4}}\\
&\le C \left( \frac{\mu(B)}{\mu(B_1)} \right)^{3/4} \mathrm{diam}(B) \left( \Barint_B |\nabla_H u|\;\mathrm{d}\mu\right)\\
&\le C \left(\frac{ \mu(B)}{\mu(B_1)}\right)^{\frac{3}{4}} \left( \frac{\mu(B_2)}{\mu(B)} \right) \mathrm{diam}(B_2) \Barint_{B_2} |\nabla_H u|\;\mathrm{d}\mu.
\end{align*}
The proof is complete.
\end{proof}

\begin{prop}[Weak Sobolev-Poincar\'{e} inequality]\label{p:weakPI}
There exists $\delta_0>0$ such that for every $\varepsilon$-net $Y_{M'}$ on $M'$ as above with $\varepsilon<\delta_0$, there are constants $0<C<\infty$ and $1\leq c<\infty$, depending only on $M'$, such that
\begin{displaymath}
\left(\Barint_{B(x,\varepsilon)}|u-u_{B(x,\varepsilon)}|^{\frac{4}{3}}\;\mathrm{d}\mu_{{M}}\right)^{\frac{3}{4}} \leq C \varepsilon \Barint_{B(x,c\varepsilon)}|\nabla_H u|\;\mathrm{d}\mu_{{M}},
\end{displaymath}
for all $x\in Y_{M'}$ and $u\in \mathcal{C}^{\infty}(M')$.
\end{prop}

\begin{proof}
For convenience, we endow $\Heis$ with the left-invariant distance  $d_{\infty}:(x,x')\mapsto \|x^{-1}\ast x'\|$, induced by the gauge function $\|(z,t)\|:=\max\{|z|,\sqrt{|t|}\}$, $(z,t)\in \mathbb{C}\times \mathbb{R}$, which is bi-Lipschitz equivalent to the standard sub-Riemannian distance $d_{\Heis}$. This gives explicit information on the shape of balls.

We fix a finite collection of sets $\{V_1,\ldots, V_n\}$ in $M$ which cover $M'$ so that each $V_i$ can be mapped $L$-bi-Lipschitzly to a domain in $(\Heis,d_\infty)$ and $V_i \cap M'$ has positive mass. Let $r_0$ be as in the proof of Proposition \ref{p:biLip_metric}, and fix $\delta_0 \le r_0$ small enough such that 
every $2L^2\varepsilon$-ball in $M$ intersected with $M'$ lies inside one of the sets $V_1,\ldots,V_n$ provided that $\varepsilon < \delta_0$.

The Darboux chart map pushes forward the Popp measure on $M'$ (induced by $g_M$) to the standard volume on $\Heis$. Let $Y_{M'}$ be a $\varepsilon$-net on $M'$ as before and let $x$ be a point in $Y_{M'}$. As previously explained there exists a neighborhood $V_i$ of $x$ which under a Darboux chart map $\Psi$ is mapped $L$-bi-Lipschitzly to a domain $U_i$ in $\Heis$. An argument similar to the one in the proof of Lemma \ref{l:incl_ball_PI} shows that it suffices to verify the inequality
\begin{align}\label{eq:PI_for_ball}
\left(\Barint_{A}|u-u_{A}|^{\frac{4}{3}}\;\mathrm{d}\mu_{\Heis}\right)^{\frac{3}{4}}\leq C \varepsilon \Barint_{2A}|\nabla_H u|\;\mathrm{d}\mu_{\Heis}
\end{align}
for all  $u\in \mathcal{C}^{\infty}(U_i)$, where
\begin{displaymath}
A=B_{d_{\infty}}(\Psi(x),L\varepsilon)\cap \Psi(M'\cap V_i), \qquad 2A= B_{d_{\infty}}(\Psi(x),2L\varepsilon)\cap \Psi(M'\cap V_i),
\end{displaymath}
and $C$ is a finite constant depending on $A$. Here we have used the fact that $d_{\infty}$-balls at $\Psi(x)$ have positive and finite $\mu_{\Heis}$-measure when intersected with $\Psi( M'\cap V_i)$.

It is not clear whether the domain $A$ is a John domain, but we will prove \eqref{eq:PI_for_ball} by applying Lemma \ref{l:incl_ball_PI} with $B_1:=A$ and $B_2:=2A$. Fix a ${\mathcal C}^{1,1}$ domain $B$ with $B_1 \subset B \subset B_2$. By \cite[Theorem 1.3]{MR2135732} $B$ is a uniform domain, and in particular a John domain in $(\Heis,d_{\Heis})$. By Theorem 1.5 and Section 6 in \cite{MR1404326} applied to $\Heis$, it follows that the desired strong $(4/3,1)$-Poincar\'{e} inequality holds for $B$. Then \eqref{eq:PI_for_ball} is a consequence of Lemma \ref{l:incl_ball_PI} and the proof is complete if we observe that there are only finitely many points in $Y_{M'}$.
\end{proof}

We are now ready to prove the local relative isoperimetric inequality for $\widetilde{M}''$.

\begin{proof}[Proof of Proposition \ref{p:localIP}]
We let $\delta_1 \leq \delta_0$, where $\delta_0$ is the parameter from Proposition \ref{p:weakPI}, which ensures that $\varepsilon$-balls centred at the points of an $\varepsilon$-net $Y_{M'}$ on $M'$ for $\varepsilon<\delta_0$ satisfy a weak $(4/3,1)$-Poincar\'{e} inequality.

We further require that $\delta_1$ is small enough such that $\pi|_{\widetilde{M}''}$ is an isometry on $c\varepsilon$-balls for $\varepsilon < \delta_1$ and $c$ as in Proposition \ref{p:weakPI}. Thus the same weak $(4/3,1)$-Poincar\'{e} inequality holds for points in the lifted net $Y$ on $\widetilde{M}''$. It then follows that there exist constants $C_I>0$ and $1<c<\infty$, depending only on $Y$ and the data of the manifold, such that for any set $E\subset \widetilde{M}''$ with locally finite perimeter,
one has
\begin{equation}\label{eq:relIP4}
C_I \left(\min \{\mu_{\widetilde{M}}(E\cap B_{\widetilde{M}''}(y,\varepsilon)),\mu_{\widetilde{M}}((\widetilde{M}''\setminus E) \cap B_{\widetilde{M}''}(y,\varepsilon))\}\right)^{\frac{3}{4}}\leq \mathcal{P}(E,B_{\widetilde{M}''}(y,c\varepsilon)),
\end{equation}
for all $y \in Y$. The proof in the sub-Riemannian setting follows the same methods as in the Euclidean case, see \cite{MR2979606} and for instance \cite[Corollary 1.29]{MR775682}: Inequality \eqref{eq:relIP4} follows by applying the Poincar\'{e} inequality to the function $u=\chi_E$. Since $\chi_E$ is not smooth but merely of bounded variation, this requires an approximation result of BV functions by smooth functions. This is given by \cite[Proposition 2.4]{MR2979606} in the setting of sub-Riemannian contact manifolds. Finally we notice that $\mathcal{P}(E,\Omega)$ is the total variation of the characteristic function $\chi_E$. The fact that we only have a \emph{weak} Poincar\'{e} inequality accounts for the enlarged ball on the right-hand side of the above isoperimetric inequality.

To conclude, we observe that for all $d>4$, we have $3/4< (d-1)/d$. Since the considered volumes are finite, we see that the desired weak local relative $d$-dimensional isoperimetric inequality follows from the above version.
\end{proof}

\subsection{Proof of the global isoperimetric inequality}
\label{sec:globalisoperimetric}

Kanai \cite[Lemma 4.5]{MR792983} established the transition from rough and local isoperimetric inequalities to global isoperimetric inequalities for {Riemannian} manifolds with lower bound on the Ricci curvature. In this section, we present a similar transition in the abstract setting of metric measure spaces. Our results apply in particular to the case of the sub-Riemannian manifold $\widetilde{M}''$.

Throughout this section we let $(X,\mu)$ be a metric measure space, which we further assume to be equipped with a {\it perimeter measure} $\mathcal P$. The perimeter measure ${\mathcal P}$ should act on pairs $E$ and $\Omega$, where $E$ is measurable and $\Omega$ is open. Further, ${\mathcal P}(\cdot,\Omega)$ should be a Borel measure for each $\Omega$, and $\Omega \mapsto {\mathcal P}(E,\Omega)$ should be monotonic with respect to set inclusion. A set $E$ is said to be of {\it finite perimeter} if ${\mathcal P}(E,X) < \infty$.

\begin{definition}
Let $(X,\mu)$ be a metric measure space equipped with a notion of perimeter $\mathcal{P}$ as above. We say that $X$ satisfies a \emph{weak relative $d$-dimensional isoperimetric inequality at scale $\varepsilon>0$ in $Y\subset X$} if there exist finite constants $c,C\geq 1$ such that
\begin{displaymath}
\left(\min\{\mu(E\cap B(x,\varepsilon)),\mu((X\setminus E)\cap B(x,\varepsilon))\}\right)^{\frac{d-1}{d}} \le C \mathcal{P}(E,B(x,c\varepsilon))
\end{displaymath}
for all $x \in Y$ and for all non-empty relative compact domains $E\subset X$ of finite perimeter. We say that $X$ satisfies a \emph{$d$-dimensional isoperimetric inequality} if there exists a constant $0<C<\infty$, such that
\begin{displaymath}
\mu(E)^{\frac{d-1}{d}} \leq C \mathcal{P}(E,X)
\end{displaymath}
for all non-empty relative compact domains $E\subset X$ of finite perimeter.
\end{definition}

Our goal in this section is to prove the following theorem.

\begin{thm}\label{t:globalIPgeneral}
Let $(X,\mu)$ be a metric measure space equipped with a perimeter function $\mathcal{P}$. Suppose that there exists $d\geq 1$ and $\varepsilon>0$ such that
\begin{enumerate}
\item\label{i:2rdcond} $X$ contains a maximal $\varepsilon$-separated net $Y$ with a $d$-dimensional isoperimetric inequality,
\item\label{i:1rdcond} $X$ satisfies a weak relative $d$-dimensional isoperimetric inequality at scales $\varepsilon$ and $3\varepsilon$ in $Y$,
\item\label{i:3rdcond} $\sup_{y \in Y} \sharp (B(y,r)\cap Y)<\infty$ for all $r>0$,
\item\label{i:4rdcond} there are constants $0<c_-$, $c_+ <\infty$ such that $c_- \leq \mu(B(y,\varepsilon)) \leq c_+$ for all $y\in Y$.
\end{enumerate}
Then $X$ satisfies a $d$-dimensional isoperimetric inequality.
\end{thm}

The proof below reveals that it suffices in fact to assume that the estimate in \eqref{i:3rdcond} holds for $r\simeq \varepsilon$, where the exact value of $r$ depends on the constant $c$ from the weak relative isoperimetric inequality. Adapting the argument given in \cite[Lemma 2.3]{MR792983}, one can see that such a weakened version of \eqref{i:3rdcond} holds true as soon as $X$ satisfies a weak condition on the volume of balls as in \eqref{i:4rdcond} -- at suitable scales depending on $\varepsilon$.


\begin{proof}[Proof of Theorem \ref{t:globalIPgeneral}]
Let $E$ be an arbitrary non-empty relatively compact domain in $X$ with finite perimeter. We wish to show that
\begin{displaymath}
\mu(E)^{\frac{d-1}{d}} \leq C \mathcal{P}(E,X)
\end{displaymath}
for a universal constant $C$ that does not depend on $E$. The strategy is to use the large-scale information provided by the rough isoperimetric inequality for the net $Y$ combined with the local information provided by the weak local relative isoperimetric inequality. Without loss of generality, we may assume that the relative isoperimetric inequality holds with the same constants at scale $\varepsilon$ and at scale $3\varepsilon$.

 We observe further that  it suffices to consider points of $Y$ that lie in the open $\varepsilon$-neighborhoood $N_{\varepsilon}(E)$ of $E$. We divide these points into two categories by setting
$$
S := \{y\in Y:\; \mu(E\cap B(y,\varepsilon))>\tfrac{1}{2} \mu(B(y,\varepsilon))\}
$$
and
$$
P_0 := \{y\in Y\cap N_{\varepsilon}(E):\; \mu(E\cap B(y,\varepsilon))\leq \tfrac{1}{2} \mu(B(y,\varepsilon))\}.
$$
We will apply the weak local relative isoperimetric inequality to points in $P_0$ and in the combinatorial boundary $\partial S$. For the set $S$ we will also apply
the rough isoperimetric inequality. Note that relative compactness of $E$ ensures that the cardinality of $S$ is finite.

By maximality of the net, the $\varepsilon$-balls centered at points in $S\cup P_0$ cover the set $E$ and thus
\begin{equation}\label{eq:vol_upper_sum}
\mu(E) \leq \sum_{y\in P_0} \mu(B(y,\varepsilon)\cap E) + \sum_{y\in S} \mu(B(y,\varepsilon)\cap E).
\end{equation}
If $y$ belongs to $P_0$, then by definition and by the weak local relative isoperimetric inequality at scale $\varepsilon$,
\begin{displaymath}
\mu(B(y,\varepsilon)\cap E)^{\frac{d-1}{d}} \leq \min\{\mu(B(y,\varepsilon)\cap E),\mu(B(y,\varepsilon)\cap (X\setminus E))\}^{\frac{d-1}{d}}\leq C \mathcal{P}(E,B(y,c\varepsilon)).
\end{displaymath}
Summing over all such points $y$, we find that
\begin{align*}
\sum_{y\in P_0} \mu(E\cap B(y,\varepsilon))&\leq \left(\sum_{y\in P_0}\left(\mu(E\cap B(y,\varepsilon))\right)^{\frac{d-1}{d}}\right)^{\frac{d}{d-1}}\\
&\leq \left(C \sum_{y\in P_0}\mathcal{P}(E,B(y,c\varepsilon)) \right)^{\frac{d}{d-1}}\\
&\leq \left(\nu C \mathcal{P}(E,X)\right)^{\frac{d}{d-1}},
\end{align*}
where the constant $\nu$ is derived from assumption \eqref{i:3rdcond}, which controls the overlap of $c\varepsilon$-balls centred in points of the net; see the argument in Remark \ref{rem:net_point_overlap}. If $S$ is empty, then this estimate combined with \eqref{eq:vol_upper_sum} gives the desired bound. Otherwise we estimate the sum over the points in $S$ as follows:
\begin{equation}\label{eq:S_sum}
\sum_{y\in S} \mu(B(y,\varepsilon)\cap E) \leq c_{+} \cdot \sharp S \leq \left(c_{+}^{\frac{d-1}{d}} C \sharp \partial S\right)^{\frac{d}{d-1}}.
\end{equation}
If $y$ is a point in $\partial S$, then by definition, $y$ does not belong to $S$, but there is a point $s\in S$ such that $d(y,s)\leq 2\varepsilon$.
Since the $\varepsilon$-ball centred at $s$ intersects $E$ significantly, but the corresponding ball centred at $y$ does not, we can show that the slightly enlarged ball $B(y,3\varepsilon)$ intersects both $E$ and its
complement in sets of large mass. Indeed, for $y$ and $s$ as above, one has
\begin{displaymath}
B(s,\varepsilon)\cap E \subseteq B(y,3\varepsilon) \cap E
\end{displaymath}
and thus
\begin{equation}\label{eq:vol_est1}
\mu(B(y,3\varepsilon) \cap E)\geq \mu(B(s,\varepsilon)\cap E)> \tfrac{1}{2}\mu(B(s,\varepsilon))\geq \tfrac{c_{-}}{2}.
\end{equation}
On the other hand, to estimate $B(y,3\varepsilon) \cap (X\setminus E)$ from below, we observe that since $y\notin S$, the point $y$ either lies outside a $\varepsilon$-neighborhood of $E$,
in which case $B(y,\varepsilon)$ is entirely contained in $X\setminus E$ and we have
\begin{equation}\label{eq:vol_est2}
\mu(B(y,3\varepsilon) \cap (X\setminus E)) \geq \mu(B(y,\varepsilon) \cap (X\setminus E))\geq \mu(B(y,\varepsilon))\geq c_{-},
\end{equation}
or $y$ belongs to $P_0$. If the latter happens, then
\begin{equation}\begin{split}\label{eq:vol_est3}
\mu(B(y,3\varepsilon) \cap (X\setminus E)) &\ge \mu(B(y,\varepsilon) \cap (X\setminus E))\\
&= \mu(B(y,\varepsilon))-\mu(B(y,\varepsilon) \cap E) \\
&\ge \tfrac{1}{2}\mu(B(y,\varepsilon))\geq \tfrac{c_{-}}{2}.
\end{split}\end{equation}
Combining \eqref{eq:vol_est1}, \eqref{eq:vol_est2} and \eqref{eq:vol_est3} with the weak local relative isoperimetric inequality at scale $3\varepsilon$, we find that
\begin{displaymath}
\left(\tfrac{c_{-}}{2}\right)^{\frac{d-1}{d}}\leq \min\{\mu(B(y,3\varepsilon)\cap E),\mu(B(y,3\varepsilon)\cap (X\setminus E))\}^{\frac{d-1}{d}}
\leq C \mathcal{P}(E,B(y,3c\varepsilon)).
\end{displaymath}
Since this estimate holds uniformly for all $y \in \partial S$, we conclude that
\begin{align*}
\sharp\partial S \leq \left(\tfrac{2}{c_{-}}\right)^{\frac{d-1}{d}} C \sum_{y\in \partial S} \mathcal{P}(E,B(y,3c\varepsilon))
\leq \left(\tfrac{2}{c_{-}}\right)^{\frac{d-1}{d}} C \nu \mathcal{P}(E,X),
\end{align*}
 where $\nu$ is again a finite constant which controls the overlap of balls, guaranteed by assumption \eqref{i:3rdcond}.

We insert this estimate in \eqref{eq:S_sum}, and return to the volume estimate \eqref{eq:vol_upper_sum} at the beginning of the proof, which now reads
\begin{displaymath}
\mu(E) \leq C \mathcal{P}(E,X)^{\frac{d}{d-1}}
\end{displaymath}
for a suitable universal constant $C$. This concludes the proof.
\end{proof}

As an application of Theorem \ref{t:globalIPgeneral}, we obtain the following statement relevant for the proof of our main result.

\begin{cor}[Global isoperimetric inequality]\label{p:isoperim_global}
The manifold $\widetilde{M}''$  satisfies a $d$-dimensional isoperimetric inequality for $d>4$, that is, there exists a constant $C>0$ such that
\begin{displaymath}
\left(\mu_{\widetilde{M}}( E)\right)^{\frac{d-1}{d}}\leq C \mathcal{P}(E,\widetilde{M}'')
\end{displaymath}
for all non-empty relatively compact domains $E\subset \widetilde{M}''$ with piecewise $\mathcal{C}^1$-boundary.
\end{cor}
Here $\mathcal{P}(E,\widetilde{M}'')$ denotes the (sub-Riemannian) perimeter of $E$ in $\widetilde{M}''$.

\begin{proof}
We verify that $X=\widetilde{M}''$ endowed with the metric $d_{\widetilde{M}''}$ and the measure $\mu_{\widetilde{M}}$ fulfills the assumptions of Theorem \ref{t:globalIPgeneral}.

We choose $\varepsilon< \min \{\varepsilon_0,\varepsilon_1, \delta_1/3\}$, where the bound for $\varepsilon$ is given by constants that have appeared earlier.
The fact that $\varepsilon< \varepsilon_0$ allows us by Lemma \ref{l:unif}  and Proposition
\ref{p:roughIP} to choose a maximal $\varepsilon$-separated net $Y$ on $\widetilde{M}''$ which satisfies condition \eqref{i:3rdcond} in Theorem \ref{t:globalIPgeneral} and fulfills a rough $d$-dimensional isoperimetric inequality for $d>4$. Since $\varepsilon$ has also been chosen smaller than $\varepsilon_1$, Lemma \ref{eq:volume_comparison} yields assumption \eqref{i:4rdcond} in Theorem \ref{t:globalIPgeneral}. Finally, the fact that $3\varepsilon$ is smaller than $\delta_1$ ensures that $\widetilde{M}''$ satisfies a weak relative $d$-dimensional isoperimetric inequality at scales $\epsilon$ and $3\varepsilon$ by Proposition \ref{p:localIP}. The claim follows.
\end{proof}

\subsection{Computation of the capacity at infinity}\label{ss:cap}
The goal of this section is to prove that the $4$-capacity in $\widetilde{M}$ of a closed ball (or a more general compact set) in $\widetilde{M}''$ is positive. We will follow a standard proof relying on the isoperimetric and coarea inequalities. We first formulate a suitable horizontal coarea formula on sub-Riemannian contact manifolds.

\begin{prop}[Coarea formula]\label{p:coarea}
Let $N$ be a contact sub-Riemannian $3$-manifold. Then, for all
$u\in \mathcal{C}^{3}(N)$, one has that
\begin{equation}\label{e:coarea}
\int_{N} |\nabla_H u| \;\mathrm{d}\mu_{N} = \int_{-\infty}^{\infty} \mathcal{P}_N(\{x\in N: u(x)> t\}, N)\,\mathrm{d}t,
\end{equation}
where $|\nabla_H u|=  \sqrt{g_{N}(\nabla_H u,\nabla_H u)}$.
\end{prop}

Here, as before, $\mu_N$ denotes the Popp volume measure on $N$ induced by the metric $g_N$, while $\mathcal{P}_N$ denotes the horizontal perimeter. Equation \eqref{e:coarea} has been established in \cite[Theorem
4.2]{MR1865002} for Lipschitz functions in Carnot-Carath\'{e}\-odory spaces, that is,
for the case when the manifold is $\mathbb{R}^n$ with a
sub-Riemannian metric. While we have to consider other manifolds as well, the coarea formula for $\mathcal{C}^3$ functions suffices for our purposes. This case is considerably easier to prove as it follows directly from the Riemannian coarea formula. The idea of using a Riemannian coarea formula to derive a statement in the sub-Riemannian setting is not new; see \cite[Theorem 7.2.2]{bMa02} (for sub-Riemannian groups).


\begin{proof}[Proof of Proposition \ref{p:coarea}]

We assume first that the contact structure of $N$ is co-orientable.
This allows to promote the sub-Riemannian metric $g_N$ to a Riemannian metric on $N$ by declaring the Reeb vector field orthonormal to the distribution $HN$. We continue to write this metric as $g_N$. Recall that the Riemannian volume associated to $g_N$ agrees with the Popp volume. Given a $\mathcal{C}^{3}(N)$ function $u$, we denote by $\nabla u$ the usual, Riemannian, gradient of $u$. Let us further abbreviate
\begin{displaymath}
E_t:= \{x\in N:\; u(x)>t\}\quad \text{and}\quad \Sigma_t:= \{x\in N:\; u(x)=t\},\quad  t\in \mathbb{R}.
\end{displaymath}
As $u$ is $\mathcal{C}^3$, it follows by Sard's theorem and the discussion in \cite[\S 2.3]{MR3412382}, \cite[(2.9)]{MR2979606} that for almost every $t\in \mathbb{R}$ (for the regular values of $u$), one has
\begin{equation}\label{eq:perim_smooth}
\mathcal{P}_N(E_t,N)= \int_{\Sigma_t} |\nu| \;\mathrm{d}\sigma^2,
\end{equation}
where $\sigma^2$ is the Riemannian measure on $\Sigma_t$ and $\nu$ the orthogonal projection to $HN$ of the unit vector field $n$ that is normal to $\Sigma_t$.  Let us fix such a regular value $t$. Since $\Sigma_t$ is a level set of $u$ and $t$ is regular, a unit vector field normal to $\Sigma_t$ is given by $n:=\nabla u/|\nabla u|$. Here $|\cdot|$ is computed with respect to $g_N$.

 Next, the Riemannian coarea formula, as stated for instance in \cite[Corollary I.3.1]{MR1849187}, says that
\begin{displaymath}
\int_N |\nabla u(x)|\phi(x)\;\mathrm{d}\mu_N= \int_{-\infty}^{\infty} \int_{\Sigma_t} \phi(y) \;\mathrm{d}\sigma^2(y)\;\mathrm{d}t,
\end{displaymath}
for all nonnegative measurable functions $\phi$ on $N$. We apply this to the function $\phi= h|\nabla_H u|/|\nabla u| $ for some nonnegative measurable function $h$ on $N$.
It follows that
\begin{equation}\label{eq:coarea_step}
\int_N h |\nabla_H u|\;\mathrm{d}\mu_N = \int_{-\infty}^{\infty} \int_{\Sigma_t} h \left|\frac{\nabla_H u}{|\nabla u|}\right| \; \mathrm{d}\sigma^2\;\mathrm{d}t.
\end{equation}
We observe that $\nabla_H u/|\nabla u|$ agrees with the orthogonal projection $\nu$ of the unit normal $n$ to $HN$. The desired coarea formula (in the case of a co-orientable contact structure) then follows from \eqref{eq:perim_smooth} and \eqref{eq:coarea_step} with $h\equiv 1$.

Since a general contact $3$-manifold can be covered by open subsets restricted to which the horizontal distribution is orientable, the general case can be proved by a partition of unity argument.
\end{proof}
%
%

\begin{definition}\label{d:cap}
Let $N$ be a contact sub-Riemannian $3$-manifold. The $p$-capacity, $1<p<\infty$, of a compact set $C\subseteq N$ is defined as
\begin{displaymath}
\mathrm{cap}_p(N,C)= \inf \int_N |\nabla_H u|^p \;\mathrm{d}\mu_N,
\end{displaymath}
where the infimum is taken over $u\in \mathcal{C}_0^{\infty}(N)$ with $u|_C\geq 1$, and the norm $|\cdot|$ is defined using the sub-Riemannian metric $g_N$.
The pair $(N,C)$ is called a \emph{condenser}.
\end{definition}

In the current section, we apply the above definition with $N=\widetilde{M}$. In Section \ref{s:morphism}, we will apply it in the case where $N$ is an open subset of $\mathbb{H}$.


\begin{prop}\label{p:4hyperbolic}
Let $K \subset \widetilde{M}''$ be a compact set of positive $\mu_{\widetilde{M}}$ measure. Then
\begin{equation}\label{e:cap-pos}
\mathrm{cap}_4(\widetilde{M}, K) >0.
\end{equation}
\end{prop}

We proceed as in \cite[Theorem 1.3]{MR2832708}, yet we work on $\widetilde{M}$ at first, and only deal with $\widetilde{M}''$ when the isoperimetric inequality comes into play, rather than deriving for instance a Sobolev inequality on $\widetilde{M}''$ in greatest possible generality. However, the reader will surely recognize in what follows arguments similar to those used in the proof of Sobolev inequalities; see for instance \cite{Pansu}.

\begin{proof}

We fix $u\in \mathcal{C}^{\infty}_0(\widetilde{M})$ such that $u|_K \geq 1$. We have to find a uniform positive lower bound for $\int_{\widetilde{M}} |\nabla_H u|^4 \,\mathrm{d}\mu_{\widetilde{M}}$. The coarea formula will naturally lead to an integral of $|\nabla_H u|$, rather than $|\nabla_H u|^4$, but this issue can be solved by applying H\"older's inequality with a suitable exponent.
For any $\gamma>1$, it holds that
\begin{equation}\label{eq:estK}
\left(\mu_{\widetilde{M}}(K) \right)^{\frac{d}{d-1}-\frac{4}{3}}\leq \left(\int_K |u|^{\gamma \frac{d}{d-1}}\,\mathrm{d}\mu_{\widetilde{M}}\right)^{\frac{d}{d-1}-\frac{4}{3}}.
\end{equation}
We will later choose the exponent $\gamma$ appropriately depending on $d$.
The classical real-variables inequality $\left( \delta^{-1} \int_0^\infty s^{1/\delta-1} F(s) \, ds \right)^\delta\leq \int_0^\infty F(s)^\delta \, ds $, valid for decreasing functions $F$ and $0<\delta\le 1$ (see, for instance, \cite[(3.34)]{MR1800917} or \cite[(II.2.2)]{MR1849187}) yields
\begin{align*}
\left(\int_{\widetilde{M}''}
\left(|u|^{\gamma}\right)^{\frac{d}{d-1}}\,\mathrm{d}\mu_{\widetilde{M}}\right)^{\frac{d-1}{d}}&\leq
\int_0^{\infty} \left(\mu_{\widetilde{M}}(x\in \widetilde{M}'':\,
|u(x)|^{\gamma}> s\right)^{\frac{d-1}{d}}\,\mathrm{d}s
\end{align*}
when applied to $F(s) = \mu_{\widetilde{M}}(\{x \in \widetilde{M''}:|u(x)|^\gamma>s\})$ and $\delta = (d-1)/d$.
Using the isoperimetric inequality, this can be further estimated from above by
\begin{align*}
C \int_0^{\infty} \mathcal{P}(\{x\in \widetilde{M}'':\, |u(x)|^{\gamma}>s\},\widetilde{M}'')\,\mathrm{d}s.
\end{align*}

If $\gamma>3$, then $|u|^{\gamma}\in \mathcal{C}^3$ since $u\in \mathcal{C}^{\infty}_0(\widetilde{M})$.
Thus we can apply the coarea formula to $N=\widetilde{M}''$ and $|u|^{\gamma}$. We obtain
\begin{equation}\label{eq:intermediate_est}
\left(\int_{\widetilde{M}''} \left(|u|^{\gamma}\right)^{\frac{d}{d-1}}\,\mathrm{d}\mu_{\widetilde{M}}\right)^{\frac{d-1}{d}}\leq C \int_{\widetilde{M}''}
|\nabla_H(|u|^{\gamma})|\,\mathrm{d}\mu_{\widetilde{M}}.
\end{equation}
For the rest of the computation, we fix
\begin{displaymath}
\gamma:= \frac{4(d-1)}{d-4}.
\end{displaymath}
Since $d>4$, it holds that $\gamma > 3$ as required. Moreover, $\gamma$ is chosen so that
\begin{displaymath}
\frac{\gamma}{\gamma-1}\frac{d}{d-1}= \frac{4}{3},
\end{displaymath}
Then H\"older's inequality and \eqref{eq:intermediate_est}
yield
\begin{displaymath}
\left(\int_{\widetilde{M}''} \left(|u|^{\gamma}\right)^{\frac{d}{d-1}}\,\mathrm{d}\mu_{\widetilde{M}}\right)^{\frac{d-1}{d}-\frac{3}{4}}\leq C \gamma \left(\int_{\widetilde{M}''} |\nabla_H u|^4\,\mathrm{d}\mu_{\widetilde{M}}\right)^4.
\end{displaymath}
Returning to \eqref{eq:estK}, we have found that
\begin{displaymath}
\left(\mu_{\widetilde{M}}(K)\right)^{\frac{d-1}{d}-\frac{3}{4}}\leq C \gamma \left(\int_{\widetilde{M}''} |\nabla_H u|^4\,\mathrm{d}\mu_{\widetilde{M}}\right)^4\leq C \gamma \left(\int_{\widetilde{M}} |\nabla_H u|^4\,\mathrm{d}\mu_{\widetilde{M}}\right)^4.
\end{displaymath}
Taking the infimum over all such $u$ completes the proof.
\end{proof}

The following notions are standard in the Riemannian setting.

\begin{definition}
We say that a contact sub-Riemannian $3$-manifold $N$ is \emph{$p$-parabolic}, $1<p<\infty$, if
\begin{displaymath}
\mathrm{cap}_p(N,C)=0
\end{displaymath}
for all compact sets $C\subseteq N$. A manifold that is not $p$-parabolic is called \emph{$p$-hyperbolic}.
\end{definition}

In this language, Proposition \ref{p:4hyperbolic} states that $\widetilde{M}$ is $4$-hyperbolic. By way of contrast, it is well known \cite[p.130]{MR1630785} that the sub-Riemannian Heisenberg group $\Heis$ is $4$-parabolic.

Next, we will introduce some machinery of nonlinear potential theory which, in coordination with the above hyperbolicity and parabolicity results, will complete the proof of Theorem \ref{thm:main}.

\subsection{Nonlinear potential theory}
\label{sec:potentialtheory}

In this section and the next, we give a brief digression into some aspects of nonlinear potential theory on sub-Riemannian manifolds.
The main goal of this section is to conclude from the
$4$-hyperbolicity of $\widetilde M$ the existence of a positive
nonconstant supersolution to the $4$-harmonic equation, and the existence
of a Green's function for the $4$-Laplacian at every point of
$\widetilde M$.

For an introduction to the classical Euclidean nonlinear potential theory, we refer the reader to \cite{MR1207810}. For a
discussion of $\mathcal{A}$-harmonic functions in the \emph{Riemannian} setting, see \cite{MR2039956}.
Nonlinear potential theory on Carnot groups has been initiated in \cite{MR1630785}. An in-depth
study of $Q$-harmonic functions on sub-Riemannian manifolds is part of \cite{CLDO}. Nonlinear potential theory in metric measure spaces of bounded geometry has been discussed in
\cite{MR2867756}. Here we will merely provide the results that are needed to prove our main theorem in the setting of manifolds modelled on the Heisenberg group.

Let $(N,HN,g_N)$ be a sub-Riemannian contact $3$-manifold. The definition of harmonic functions $u:N \to \mathbb{R}$ requires the notion of a horizontal gradient, divergence, and a Laplacian on $N$.
 The divergence of a smooth vector field has been defined in Definition \ref{eq:smoothDiv}, the horizontal gradient in Definition \ref{d:horiz_grad}. We now extend these notions to the nonsmooth case.

We say that a horizontal vector field $\nabla_H u \in
L^1_{loc}(N)$ is a \emph{weak horizontal gradient} of $u\in
L_{loc}^1(N)$ if
\begin{displaymath}
\int_N g_N(\nabla_H u, \Phi)\mathrm{d}\mu_N = - \int_N u\,
\mathrm{div}_N\Phi \;\mathrm{d}\mu_N
\end{displaymath}
for all smooth compactly supported horizontal vector fields $\Phi$ on $N$.

\begin{remark}
If a continuous map $f:\Heis \to N$ has weak horizontal
derivatives in $L^4_{loc}$, then $u\circ f\in
HW_{loc}^{1,4}(\Heis)$ for every smooth function $u:N\to
\mathbb{R}$, and the weak horizontal gradient of $u\circ f$
exists and agrees with $X(u\circ f)X + Y(u\circ f)Y$.
\end{remark}

We say that $\mathrm{div}_N V \in L_{loc}^1(N)$ is  a \emph{weak
divergence} of a locally integrable horizontal vector field $V$ on
$N$ if
\begin{displaymath}
\int_N g_N(\nabla_H \varphi, V)\mathrm{d}\mu_N = - \int_N \varphi\,
\mathrm{div}_N V \;\mathrm{d}\mu_N
\end{displaymath}
for all $\varphi\in \mathcal{C}_0^{\infty}(N)$.

\begin{definition}\label{d:lap}
A continuous $HW_{loc}^{1,4}$-function $u:N \to \mathbb{R}$
is said to be
\emph{$4$-harmonic}, if the equation
\begin{displaymath}
\mathrm{div}_N\left(g_N(\nabla_H u, \nabla_H u)
\nabla_H u\right)=0
\end{displaymath}
holds in a weak sense.
\end{definition}

In the following, we will also use a generalization of this
concept. We consider operators $\mathcal{A}: HN \to HN$ for which there exist
constants  $0<\alpha \leq \beta <\infty$ such that
\begin{enumerate}
\item $\mathcal{A}_x:H_xN \to H_x N$ is continuous for almost
every $x$, \item $x \mapsto \mathcal{A}_x(V)$ is measurable for
all horizontal measurable vector fields $V$, \item for almost every $x\in
N$ and all $h\in H_xN$:
\begin{enumerate}
\item\label{i:a} $g_N(\mathcal{A}_x(h),h)\geq \alpha g_N(h,h)^{2}$,
\item\label{i:b} $g_N(\mathcal{A}_x(h),\mathcal{A}_x(h))^{1/2}\leq \beta g_N(h,h)^{3/2}$,
\item $g_N(\mathcal{A}_x(h_1)-\mathcal{A}_x(h_2),h_1-h_2)>0$ for $h_1\neq h_2$,
\item $\mathcal{A}_x(\lambda h)=|\lambda|^{2}\lambda \mathcal{A}_x(h)$ for all $\lambda \in \mathbb{R}\setminus \{0\}$.
\end{enumerate}
\end{enumerate}
Here and in what follows we have written ${\mathcal{A}}_x(h) := \mathcal{A}(x,h)$ for $h \in H_xN$.
We will call such $\mathcal A$ {\it operators of type $4$} on $N$.

\begin{definition}\label{d:harmonic_sol}
A $HW_{loc}^{1,4}$-function $u:N \to \mathbb{R}$  is called \emph{solution of the
$\mathcal{A}$-harmonic equation}, or for short an {\it ${\mathcal A}$-solution} if
\begin{equation*}
-\mathrm{div}_N\left(\mathcal{A}( \nabla_H u)\right)=0
\end{equation*}
holds in the weak sense for a
suitable
$\mathcal{A}:HN \to HN$ as above. ${\mathcal A}$-\emph{sub-} and
\emph{super-solutions} can be defined accordingly using the signs
$\leq$ and $\geq$, and nonnegative test functions. If the solution $u$ is continuous, it is called
$\mathcal{A}$-harmonic.
\end{definition}

The {\it standard operator of type $4$} is ${\mathcal A}_x(h) = g_N(h,h)h$;
continuous solutions to the $\mathcal A$-harmonic equation for this $\mathcal A$ are precisely the $4$-harmonic functions of Definition \ref{d:lap}.

We will also encounter solutions with a singularity. Let $\Omega$
be a relatively compact domain in $N$ and $y$ a point in $N$. We
say that a positive function
\begin{displaymath}
G=G(\cdot,y)\in \mathcal{C}(\Omega \setminus \{y\})\cap
HW_{loc}^{1,4}(\Omega \setminus \{y\})
\end{displaymath}
is a \emph{Green's function in $\Omega$ with pole $y$ for the
$\mathcal{A}$-harmonic equation} if
\begin{enumerate}
\item $\lim_{x\to z}G(x)=0$ for all $z\in \partial \Omega$, \item
$\int_{\Omega} g_N(\nabla_H \varphi,\mathcal{A}(\nabla_H
G))\;\mathrm{d}\mu_N = \varphi(y)$ for all $\varphi \in
\mathcal{C}_0^{\infty}(\Omega)$.
\end{enumerate}
Note that such $G(\cdot,y)$ is $\mathcal{A}$-harmonic in $\Omega
\setminus \{y\}$.

We say that a function $G=G(\cdot,y)$ is a \emph{Green's function
in $N$ with pole $y$ for the $\mathcal{A}$-harmonic equation} if
there exists an exhaustion of $N$ by relatively compact domains
$\Omega_i \subset \Omega_{i+1}$, $i\in \mathbb{N}$, and associated
Green's functions $G_i$ with pole at $y$, such that $G$ equals
$\lim_{i\to \infty}G_i$ and is not identically equal to infinity.

The $4$-parabolicity of the Heisenberg group implies a
Liouville-type consequence for ${\mathcal A}$-supersolutions. The
following theorem is stated in the context of general Carnot
groups in \cite[p.\ 131]{MR1630785}, see also \cite[\S
3.2]{MR1878317} and \cite[Theorem 4]{MR1672629}.

\begin{thm}\label{t:parab_equiv0}
Let $\mathcal A$ be an operator of type $4$ on $\Heis$. Then every
nonnegative
$\mathcal A$-supersolution $u$ with
\begin{equation}\label{eq:lsc}
u(x)=\mathrm{ess}\,\liminf_{y\to x}u(y),\quad x\in \Heis,
\end{equation}
is constant.
\end{thm}

The statements in \cite{MR1630785} and \cite{MR1672629}
are formulated for so called ``superharmonic functions'' rather
than for ``supersolutions''. Yet it is not difficult to see that a
nonnegative supersolution with the property \eqref{eq:lsc} is
lower semicontinuous and fulfills the comparison principle
required in the usual definition of superharmonic functions. The
proof in the Euclidean setting, \cite[Theorem 7.16]{MR1207810},
can be easily adapted to the  Heisenberg group; the only point
worth observing is that a horizontal Sobolev function with almost
everywhere vanishing horizontal gradient has to be constant almost
everywhere. (See also \cite[Proposition 9.4]{MR2867756} for
superharmonic functions in an abstract metric measure space
setting.)

Complementing Theorem \ref{t:parab_equiv0}, we have the following
result.

\begin{thm}\label{t:parab_equiv}
If a contact sub-Riemannian $3$-manifold $N$ is $4$-hyperbolic,
then it admits a nonconstant nonnegative supersolution of the
standard operator $\mathcal A$ of type $4$, with the property that
\begin{equation}\label{eq:lsc_reg_onN}
u_N(x) = \mathrm{ess}\,\liminf_{y\to x}u_N(y),\quad x\in N.
\end{equation}
Moreover, $N$  supports a positive Green's function $G(\cdot,y)$
for $\mathcal A$ at any $y\in N$.
\end{thm}


We sketch an argument for Theorem \ref{t:parab_equiv}. The
approach is standard, and we refer to \cite[Theorem
9.22]{MR1207810} (for Euclidean spaces), \cite[Theorem
3.27]{HolDiss} (Riemannian setting), \cite[Section 3.2]{MR1878317}
and \cite[Theorems 3 and 4]{MR1672629} (sub-Riemannian setting),
\cite[Theorem 3.14]{MR1940332} (for abstract metric measure
spaces).

\begin{proof}
Since $N$ is $4$-hyperbolic, it contains a compact set $K$ whose
$4$-capacity at infinity is positive and bounded. That is, there
is a $C>0$ such that for any open set $\Omega$ with $K\subset
\Omega \Subset N$, one has
\begin{equation}\label{eq:capacityatinfinity}
\inf \int_{\Omega} \norm{\nabla_H u}^4 \;\mathrm{d}\mu_N\ge C,
\end{equation}
where the infimum is taken over $u\in C^\infty_0(\Omega)$
satisfying $u\vert_K\geq 1$. Without loss of generality, we may
restrict to \emph{nonnegative} functions $u$.

For a fixed $\Omega$, consider a minimizing sequence $u_j$ for
\eqref{eq:capacityatinfinity}. Each $u_j$ is in
$C_0^\infty(\Omega)$, and therefore in $HW^{1,4}_0(\Omega)$. One
shows that the sequence $u_j$ is a Cauchy sequence in this Sobolev
space and that the limit \emph{potential} function
$u_\Omega:=\lim_{j\rightarrow \infty} u_j \in HW^{1,4}_0$ is
nonnegative and $4$-harmonic outside of $K$. Furthermore, one
shows that
\begin{equation}
\inf \int_{\Omega} \norm{\nabla_H u}^4\;\mathrm{d}\mu_N=
\int_{\Omega} \norm{\nabla_H u_\Omega}^4\;\mathrm{d}\mu_N.
\end{equation}
Consider now an exhaustion of $N$ by domains $\Omega_j \Subset N$,
with associated potential functions $u_{\Omega_j}$. Again, one
shows that these converge
in the Sobolev space to a potential
$u_N$, now defined on all of $N$.  The limiting function is
nonnegative, and satisfies
\begin{equation}\inf \int_{N} \norm{\nabla_H u}^4\;\mathrm{d}\mu_N= \int_{N} \norm{\nabla_H u_N}^4\;\mathrm{d}\mu_N\end{equation}
In particular, $u_N$ is nonconstant. By Theorem 8.22 in
\cite{MR2867756}, we may choose an $HW_{loc}^{1,4}$ representative
of $u_N$ (which we continue to denote by the same letter) for
which \eqref{eq:lsc_reg_onN} holds. To show that $u_N$ is a
supersolution, one considers the variational kernel
\begin{displaymath}
F(x,\xi):= |\xi|^4 = g_N(\xi,\xi)^2,\quad \xi\in H_x N
\end{displaymath}
and the associated variational integral
\begin{displaymath}
I_F(u):= \int_N F(x,\nabla_H u)\;\mathrm{d}\mu_N.
\end{displaymath}
By construction, $u=u_N$ is a superminimizer for $I_F$ in the
sense of \cite[Definition 7.7]{MR2867756} and thus one shows
analogously as in the Euclidean case (\cite[Theorem
5.13]{MR1207810}) that
\begin{displaymath}
\int_N g_N\left(g_N(\nabla_H u_N,\nabla_H u_N)\nabla_H u_N,
\nabla_H v -\nabla_H u_N\right)\;\mathrm{d}\mu_N \geq 0
\end{displaymath}
for all admissible $v=u_N + \varepsilon \varphi$. This shows that
$u_N$ is a supersolution of the $4$-Laplacian on all of $N$.

To construct a Green's function, one takes a sequence of balls
$K_j = B(y, r_j)$ with $r_j \rightarrow 0$ and shows that the
global potential functions associated to $K_j$ converge, up to
renormalization, to a Green's function.
\end{proof}

\subsection{Morphism property}\label{s:morphism}

In this section, we show that if $u: N \rightarrow \R$ is a
$4$-harmonic function and $f: \Heis \rightarrow N$ is a
quasiregular mapping, then the composition $f\circ u$ is $\mathcal
A$-harmonic for a suitable operator $\mathcal A$ of type $4$ on
$\Heis$.

This so-called {\it morphism property} has been proved in
\cite[Theorem 3.14]{MR1630785} (under an additional smoothness
assumption on the mapping) and in \cite{MR1721676} (without such
assumption) for arbitrary quasiregular maps between domains in the
sub-Riemannian Heisenberg group. A morphism property for
\emph{$1$-quasiconformal} maps between equiregular sub-Riemannian
manifolds has recently been proved in \cite{CLDO}. None of these
results covers exactly the case we are interested in, on the other
hand, unlike in the setting of the mentioned results, we can rely
on an already well established theory of quasiregular mappings in
the Heisenberg group.

The {\it pullback} of an operator $\mathcal A$ of type $4$ under a
quasiregular mapping $f$ is the operator $f^\#\mathcal A$ whose
value at a point $x$ on a horizontal tangent vector $h$ is
$$
f^\# \mathcal A_x(h) := \det(D_Hf(x))^2 D_Hf(x)^{-1} \mathcal
A_{f(x)}((D_Hf(x)^{-1})^T h)
$$
if $\det D_H f(x) \neq  0$, and $f^\# \mathcal A_x$ equal to the
standard operator of type $4$ at $x$ otherwise. See Section
\ref{ss:JacQR} for the definition of the formal Jacobian $\det D_H
f$.

\begin{lemma}\label{lem:morphism1}
Let $f: V \rightarrow V'$ be a quasiregular mapping, for $V\subset \Heis$ and $V' \subset N$ domains and $N$ a sub-Riemannian $3$-manifold.
Let $\mathcal A$ be the standard operator of type $4$ in $V'$. Then $f^\#\mathcal A$ is an operator of type $4$ in $V$.
\end{lemma}

\begin{proof}
The proof goes analogously to the Euclidean
case (see \cite[Lemma 14.38]{MR1207810}), using the
characterization of quasiregularity provided in Proposition
\ref{p:qr_formal_jac} and the fact that for a nonconstant
quasiregular map $f$, a set $A$ has measure zero if and only if
$f(A)$ has measure zero (see Remark \ref{r:Lusin}).
\end{proof}

\begin{lemma}\label{lem:morphism2} Assume that $V$ and $V'$ are domains in $\Heis$, and $V''$ is a domain in an arbitrary sub-Riemannian contact $3$-manifold.
For quasiregular mappings $h:V \to V'$ and $f:V' \to V''$ and an operator $\mathcal A$ of type $4$ in $V''$ we have
$$h^\# f^\# \mathcal A = (f\circ h)^\#\mathcal A.$$
\end{lemma}

\begin{proof}
The statement is a simple computation based on the chain rule
\begin{displaymath}
D_H (f \circ h)(p) = D_H f(h(p)) \circ D_H h(p),\quad \text{for almost every }p\in V.
\end{displaymath}
The latter follows from Proposition \ref{p:chainQRHeis} applied to
the map $g$ and the components $u=f_i$, $i\in \{1,2\}$, of $f$ in
coordinates.
\end{proof}

\begin{prop}\label{prop:easyMorphism}
Let $N$ be a smooth sub-Riemannian contact $3$-manifold, $V\subset
\Heis$ and $V'\subset N$ be domains, and $\phi:V\to V'$  a
quasiconformal diffeomorphism. If $u$ is a $4$-harmonic function,
then $v:= u \circ \phi$ is $\phi^\#\mathcal A$-harmonic for the
standard operator $\mathcal A$ of type $4$. An analogous statement
holds for supersolutions.
\end{prop}

\begin{proof}
By Lemma \ref{lem:morphism1}, $\phi^\# \mathcal A$ is an operator
of type $4$. We would like to show that $v$ is a weak solution for
the $\phi^\# \mathcal A$-Laplacian. That is, for any $\Psi \in
C_0^\infty(V)$,
$$\int_V g_{\Heis}( \phi^\#\mathcal A_x(\nabla_H v), \nabla_H \Psi ) \;\mathrm{d}\mu_{\Heis}= 0.$$
We push $\Psi$ forward via $\phi$, obtaining $\phi^\#\Psi = \Psi \circ \phi^{-1}:V'\to\R$, and compute with the help of Proposition \ref{p:Jac2} that
\begin{align*}
&\int_V g_{\Heis} (\phi^\#\mathcal A_x(\nabla_H v), \nabla_H \Psi )  \;\mathrm{d}\mu_{\Heis}\\
 &= \int_V g_{\Heis} ( \det(D_H\phi(x))^2 D_H\phi(x)^{-1} \mathcal A_{\phi(x)}((D_H\phi(x)^{-1})^T \nabla_H v), \nabla_H \Psi )  \;\mathrm{d}\mu_{\Heis}\\
 &= \int_V \det(D_H\phi(x))^2  g_{\Heis} (  \mathcal A_{\phi(x)}((D_H\phi(x)^{-1})^T  \nabla_H v), (D_H\phi(x)^{-1})^T\nabla_H \Psi )  \;\mathrm{d}\mu_{\Heis}\\
 &= \int_{V'}  g_N(  \mathcal A_{x}(\nabla_H u), \nabla_H (\phi^\#\Psi))\;\mathrm{d}\mu_{N},
\end{align*}
which is equal to zero since $\phi^\# \Psi \in C^\infty_0(V')$.
\end{proof}

In the next step, we will pull $\mathcal A$ back by a quasiregular function that need not be a diffeomorphism.

\begin{prop}\label{prop:secondMorphism}
Let $U, V\subset \Heis$ be domains and $h: U \rightarrow V$ a
quasiregular mapping. If $v: V \rightarrow \R$ is $\mathcal
A$-harmonic for some operator $\mathcal A$ of type $4$, then $w =
v\circ h$ is $h^\#\mathcal A$-harmonic. An analogous statement
holds for supersolutions.
\end{prop}

This result has also been stated in \cite[Theorem 9]{MR1672629}
for $\mathcal{A}$-harmonic functions. The main technical difficulty in the proof is to push forward a test function $\Psi$ under a quasiregular mapping $h$. If $h$ was a homeomorphism, such a push-forward could be simply defined as $\Psi \circ h^{-1}$. If $h$ is not injective, it is still possible to define a function $h^\#\Psi$ which plays the role of a push-forward, but it is more difficult to verify the necessary regularity properties.
To do so, we use some terminology from
topology, for which we refer to \cite{MR1238941} or
\cite[14.9]{MR1207810}. For the moment, let us just recall that a
relatively compact domain is called a \emph{normal domain} for a
map $h$ if $h(\partial D)=
\partial h(D)$. Recall further that nonconstant quasiregular mappings on
$\Heis$ are discrete and open \cite{Da}, and by the latter
property we have $\partial h(D) \subset h(\partial D)$ for every
domain. We employ some terminology related to path
lifting. For an interval $[a,b]$ and $c\in (a,b]$, we write $I_c
= [a,c)$ if $c<b$, and $I_c=[a,b]$ if $c=b$. Given a path
$\beta:[a,b]\to \Heis$, we say that a path $\alpha: I_c \to \Heis$
is an $h$-\emph{lifting} of $\beta$ starting at a point $x\in \Heis$  if $\alpha(a)=x$ and $h\circ \alpha=\beta|_{I_c}$. We
call $\alpha$ a \emph{total $h$-lifting} if $I_c=[a,b]$.

\begin{proof}[Proof of Proposition \ref{prop:secondMorphism}] Without loss of generality, we may assume that $h$ is nonconstant and that $U$ is a normal neighborhood whose $h$-image is a ball $V$.
As in Proposition \ref{prop:easyMorphism}, we need to push forward
a test function $\Psi\in C^\infty_0(U)$. By a result of Dairbekov
\cite{MR1778673}, quasiregular mappings between domains in the
Heisenberg group are discrete and open, hence index theory is
applicable and the branch set and its image both have measure
zero. We can then define the push-forward of $\Psi$ as
$$h^\# \Psi (y) = \sum_{x\in h^{-1}y} \operatorname{index}(h, x) \Psi(x),$$
where $\operatorname{index}(h,x)$ is the local topological index of $h$ at $x$.

As in \cite[Lemma 14.30]{MR1207810}, one verifies that $h^\# \Psi\in \mathcal{C}_0(V)$ and the support of $h^\# \Psi$ is contained in $h(\mathrm{spt}\Psi)$.
While $h^\#\Psi$ is not necessarily smooth, one can show that it
is in $HW^{1, 4}_{0}(V)$. This can be done along the lines of Lemma
14.31 in \cite{MR1207810}: one works locally and proves absolute
continuity of $h^\# \Psi$ along almost every horizontal line in
$V$. These lines are first lifted under $h$  to curves in $U$, and then one
shows that almost every  such $h$-lifting is absolutely continuous.
This is the content of Lemma \ref{l:ac_lift} below.
Once the $ACL$-property of $h^\# \Psi$ is established, one shows by the same argument as in \cite[Lemma 14.31]{MR1207810}
that
\begin{displaymath}
\int_V |\nabla_H (h^\# \Psi)|^4 \;\mathrm{d}\mu_{\mathbb{H}} <\infty.
\end{displaymath}
We can then carry out the calculations as in Proposition
\ref{prop:easyMorphism} and use the density of smooth functions in
$HW^{1,4}_{loc}$ to conclude
$$\int_{V}  \langle  \mathcal A_{x}(\nabla_H v), \nabla_H h^\#\Psi \rangle\;\mathrm{d}\mu_{\Heis} =0$$ as desired.
\end{proof}

In the preceding proof we made use of the following result.

\begin{lemma}\label{l:ac_lift}
Let $h:U\to\mathbb{H}$ be a nonconstant quasiregular mapping in a domain $U\subseteq\mathbb{H}$ and let $\Psi\in\mathcal{C}^{\infty}_0(U)$. Then $h^\#\Psi\in \mathrm{ACL}(h(U))$.
\end{lemma}

Lemma \ref{l:ac_lift} can be found in \cite[Lemma 7]{MR2679532} in a more general setting. Here we specialize to the Heisenberg group and the class of quasiregular mappings.
The proof in \cite{MR2679532} is based on a series of results in other papers, which we will list below. The method of proof differs from the argument in Euclidean spaces, where it is used that quasiregular mappings have bounded inverse metric dilatation. Instead, the proof in \cite{MR2679532} makes use of a capacity estimate, which we think deserves to be better known. In \cite[Lemma 5]{MR1654140}, the following was proved (in greater generality): if $E\subset \mathbb{H}$ is connected and $G\subset\mathbb{H}$ is an open set contained in the metric $c_0 \mathrm{diam}E$-neighborhood of $E$ for a given universal constant $c_0>0$, then
\begin{equation}\label{eq:cap_lower_bound}
\left(\mathrm{cap}_4(G,E)\right)^3 \geq  c \frac{(\mathrm{diam}E)^4}{\mu_{\mathbb{H}}(G)}
\end{equation}
for an absolute constant $0<c<\infty$. (Note that the smoothness assumption on the admissible functions in our definition of capacity can be relaxed by an approximation argument, so as to make it agree with the definition given in \cite{MR1654140}.)

For the Euclidean antecedent of \eqref{eq:cap_lower_bound}, see \cite[Lemma 5.9]{MR0259114}. The estimate \eqref{eq:cap_lower_bound} is useful when coupled with a distortion inequality for quasiregular mappings and condensers.
It is straightforward to verify, see for instance \cite[Proposition 2]{MR2679532},  that for every quasiregular mapping $h:U \to \mathbb{H}$, $U\subseteq \mathbb{H}$, there exists a constant $1\leq K<\infty$, such that for every normal domain $A \subset U$ and every condenser $(A,C)$, one has
\begin{equation}\label{eq:cap_est}
\mathrm{cap}_4(A,C)\leq K N(h,A)\mathrm{cap}_4(h(A),h(C)),
\end{equation}
where $N(h,A):= \sup_{x\in \mathbb{H}}\sharp(h^{-1}(x)\cap A)$.

For the benefit of the reader we will work out in detail that part of the proof of Lemma \ref{l:ac_lift} which concerns the application of \eqref{eq:cap_lower_bound}.
We sketch the remaining part of the argument and refer the reader to the cited references for more details.

\begin{proof}[Proof of Lemma \ref{l:ac_lift}]
Throughout the proof we assume that $\mathbb{H}$ is endowed with the Kor\'{a}nyi distance
\begin{displaymath}
d(p,q):= \|p^{-1}\ast q\|_K,\quad \|(x,y,t)\|_K=\sqrt[4]{(x^2+y^2)^2 + t^2},
\end{displaymath}
which is bi-Lipschitz equivalent to the sub-Riemannian distance $d_{\mathbb{H}}$.
Let $x_0 \in \mathrm{supp} (h^\#\Psi)$ and $h^{-1}(x_0)\cap \mathrm{supp}\Psi=\{q_1,\ldots,q_s\}$. Without loss of generality, we may assume that $x_0$ is the origin. One chooses small enough normal neighborhoods $U_k:=U(q_k,h,r_1)$ around $q_k$ with $h(U_k)=B(x_0,r_1)$ as described in  \cite{MR2679532}. We may assume that all the $U_k$ are compact; cf.\ \cite[I, Lemma 4.9]{MR2210111}.

Following \cite{MR2679532}, we construct
a ``cube'' $Q$ inside $B(x_0,r_1)$ which is fibered by segments $\beta_z$ along
the flow lines of a left invariant horizontal vector field $V$, where $z$ ranges in a domain of a hyperplane transversal to $V$. The first task is to show for almost every $z\in S$ that every total $h$-lifting
 $\alpha:[a,b]\to U_k$ of the horizontal curve $\beta_z:[a,b] \to Q$  is absolutely continuous. To do so, one introduces the set function
\begin{displaymath}
\Phi(V):= \mu_{\mathbb{H}}\left(\bigcup_{k=1}^s U_k \cap h^{-1}(V\cap Q)\right),\quad V\text{ Borel}.
\end{displaymath}
In \cite[Proposition 1]{MR1375429}, it was shown that the upper volume  derivative
\begin{displaymath}
\overline{\Phi}'(z):= \limsup_{r\to  0+}\frac{\Phi(N_r(\beta_z)\cap Q)}{r^3}
\end{displaymath}
exists and is finite for almost all $z\in S$. Here, $N_r(\beta_z)$, $r>0$, denotes the metric $r$-neighborhood of $\beta_z$. We fix a point $z\in S$ where $\Phi$ has finite upper volume derivative and argue that all the $h$-liftings $\alpha$ of $\beta_z:[a,b]\to Q$ in $U_k$, $1\leq k\leq s$, are absolutely continuous.

To do this, we choose disjoint closed arcs $I_1,\ldots,I_l$ with respective lengths $\triangle_1,\ldots,\triangle_l$ on $\beta_z$ such that
\begin{displaymath}
\sum_{i=1}^l \triangle_i <\delta.
\end{displaymath}
We define $[a_i,b_i]:=\beta_z^{-1}(I_i)\subset [a,b]$. Then $E_i:= \alpha([a_i,b_i])$ is a connected set in $U_k$. We will show that $\sum_{i=1}^l d(\alpha(a_i),\alpha(b_i))\leq \sum_{i=1}^l \mathrm{diam}(E_i)$ can be made smaller than any given constant if $\delta>0$ is chosen small enough. To achieve this goal, we will construct small open neighborhoods $G_i$ of $E_i$ such that -- among other assumptions -- the  conditions for the estimate \eqref{eq:cap_lower_bound} are satisfied for the condensers $(G_i,E_i)$ for all $1\leq i\leq l$.


First, by continuity of $h$, there exists $0<r_2<c_0 \min \{\mathrm{diam}(E_i):\;1\leq i\leq 1\}$ such that $h(N_{r_2}(E_i))$ is compactly contained in $B(x_0,r_1)$. We may further assume that $r_2$ is small enough so that all the sets $h(N_{r_2}(E_i))$, $1\leq i\leq l$, are disjoint.

Second, since $h$ is an open mapping, it follows that $I_i$ is at positive distance from the boundary of $h(N_{r_2}(E_i))$ for all $1\leq i\leq l$, and so we can choose $0<r_3<\min\{\mathrm{dist}(I_i,\partial h(N_{r_2}(E_i))):\; 1\leq i\leq l\}$ such that
\begin{displaymath}
N_{r_3}(I_i)\subset h(N_{r_2}(E_i)),\quad 1\leq i\leq l.
\end{displaymath}
We may assume that $r_3<\delta$.

Third, since $U_k$ is a normal domain, for every $r<r_3$, the components of $h^{-1}(N_r(I_i))\cap U_k$ are mapped by $h$ onto $N_r(I_i)$; see \cite[I, Lemma 4.8]{MR2210111}. This shows that the $E_i$-component of $h^{-1}(N_r(I_i))\cap U_k$ is contained entirely inside $N_{c_0 \mathrm{diam}(E_i)}(E_i)$. For if this was not the case, then part of the boundary of $N_{r_2}(E_i)$ would have to be mapped inside $N_r(I_i)$, which is impossible by the choice of $r_3$.

Finally, we define $G_i$ to be the $E_i$-component of $h^{-1}(N_{r}(I_i))\cap U_k$. Note that $G_i$ is a normal domain by \cite[I, Lemma 4.7]{MR2210111} and $(G_i,E_i)$ is a condenser which satisfies the conditions for the capacity lower bound \eqref{eq:cap_lower_bound}.

The image $(h(G_i),h(E_i))$, $h(E_i)=I_i$ is again a condenser. By \eqref{eq:cap_est}, it follows that
\begin{displaymath}
c^{\frac{1}{3}}\frac{\mathrm{diam}(E_i)^{\frac{4}{3}}}{\mu_{\mathbb{H}}(G_i)^{\frac{1}{3}}}\leq \mathrm{cap}_4(G_i,E_i)\leq K N(h,G_i) \mathrm{cap}_4(h(G_i),h(E_i)).
\end{displaymath}
This implies that
\begin{align}\label{eq:arg_acl}
\notag \mathrm{diam}(E_i) &\leq c^{-\frac{1}{4}} K^{\frac{3}{4}} N(h,G_i)^{\frac{3}{4}} \cdot \left(\frac{\mu_{\mathbb{H}}(G_i)}{r^3}\right)^{\frac{1}{4}} \cdot \left(r \, \mathrm{cap}_4(h(G_i),h(E_i))\right)^{\frac{3}{4}}\\
&\leq c' \cdot \left(\frac{\Phi(N_r(\beta_z)\cap Q)}{r^3}\right)^{\frac{1}{4}} \cdot \left(r \, \mathrm{cap}_4(h(G_i),h(E_i))\right)^{\frac{3}{4}}
\end{align}
for $c'= c^{-\frac{1}{4}} K^{\frac{3}{4}} N(h,U_k)^{\frac{3}{4}}$. The proof is nearly complete if we find a constant $c''$ (which is allowed to depend on $z$, $k$ and $h$) such that
\begin{equation}\label{eq:goal_cap_upperbound}
\mathrm{cap}_4(h(G_i),h(E_i)) \leq \frac{c'' \triangle_i}{r},\quad 1\leq i\leq l.
\end{equation}
By construction, $h(G_i)=N_r(E_i)$.
 This shows that
\begin{displaymath}
w:\mathbb{H} \to \mathbb{R},\quad w(q):= \left\{\begin{array}{ll}\frac{\mathrm{dist}(q,\partial N_{r/2}(E_i)}{r/2},&w\in N_{r/2}(\beta_z)\\ 0&w\in \mathbb{H}\setminus N_r(\beta_z) \end{array} \right.
\end{displaymath}
is admissible for  $\mathrm{cap}_4(h(G_i),h(E_i))$. Indeed, by definition of the metric neighborhood,  it follows that $w(q)\geq 1$ for $q\in E_i$; and $w$ vanishes in  a neighborhood of $\partial N_r(E_i)$.
By the $1$-Lipschitz continuity of $ \mathrm{dist}(\cdot, C)$,
$$
|\nabla_H \mathrm{dist}(\cdot, C)|\leq 1
$$
for any compact set $C$. Thus
\begin{displaymath}
\mathrm{cap}_4(h(G_i),h(E_i))\leq \int_{N_r(E_i)} |\nabla_H w(q)|^4 \;\mathrm{d}\mu_{\mathbb{H}}(q) \leq 2^4 \frac{\mu_{\mathbb{H}}(N_r(E_i))}{r^4}\leq c'' \frac{\triangle_i r^3}{r^4},
\end{displaymath}
for some constant $c''$ which does not depend on $\delta$ and the choice of $I_1,\ldots,I_l$. This yields \eqref{eq:goal_cap_upperbound}. The proof of the absolute continuity of $\alpha$ concludes by \eqref{eq:arg_acl} as in  \cite{MR2679532}. It is then straightforward to verify that $h^\#\Psi$ is in $ACL$, see the arguments in \cite{MR1207810} or \cite{MR2679532}.
\end{proof}

Combining the results of this subsection with the Darboux Theorem and using Lemma \ref{lem:morphism2}, we deduce the following theorem.

\begin{thm}\label{t:morphism}
Let $f:\Heis \to   N$ be a quasiregular map from the Heisenberg
group to a smooth sub-Riemannian contact $3$-manifold $N$. If $u$
is a  $4$-harmonic function, then $w:= u \circ f$ is $f^\#
\mathcal{A}$-harmonic, where $\mathcal A$ is the standard operator
of type $4$. An analogous statement holds for supersolutions.
\end{thm}

\begin{proof}
We fix $x \in \Heis$ and a Darboux chart $\phi$ mapping an open
set of $\Heis$ to an open neighborhood of $f(x)$ in $N$. We then
write locally $u \circ f = v \circ h$ where $v = u \circ \phi$ and
$h = \phi^{-1} \circ f$. By Proposition \ref{prop:easyMorphism},
$v$ is $\phi^\#\mathcal A$-harmonic. By Proposition
\ref{prop:secondMorphism}, $w = v\circ h$ is $(\phi^{-1}\circ
f)^\# \phi^\# \mathcal A$-harmonic. Finally, Lemma
\ref{lem:morphism2} implies that $(\phi^{-1}\circ f)^\# \phi^\#
\mathcal A =
(\phi \circ \phi^{-1} \circ f)^\#\mathcal A= f^\#\mathcal A$. The proof is complete.
\end{proof}

\subsection{Application to quasiregular mappings}
\label{sec:victory}
Applying the following theorem with $N=\widetilde{M}$ concludes the proof of Theorem \ref{thm:main} by the argument given in the introduction.

\begin{thm}\label{t:no_qr}
If $f:\Heis \to N$ is quasiregular and $N$ is $4$-hyperbolic, then $f$ is constant.
\end{thm}

\begin{proof}
The proof rests on the fact that, unlike $N$, the Heisenberg group is $4$-parabolic.
If $f$ is surjective, we arrive at a contradiction by the morphism property (Theorem \ref{t:morphism}) and Theorems \ref{t:parab_equiv0} and \ref{t:parab_equiv}.

If $f$ is not constant, but $f(\Heis)$ is a strict subset of $N$,
then we choose a point $y$ on the boundary of $f(\Heis)$. The
$4$-hyperbolicity of $N$ allows us to select a  positive Green's
function $G=G(\cdot,y)$ for the $4$-Laplacian on $N$ (see Theorem
\ref{t:parab_equiv}). By Theorem \ref{t:morphism}, $G\circ f$
would be a positive nonconstant solution to an $\mathcal
A$-harmonic equation for some operator $\mathcal A$ of type $4$ on
$\Heis$. However, such solutions cannot exist (see Theorem
\ref{t:parab_equiv0}).
\end{proof}

\appendix

\section{Calculus for horizontal derivatives}\label{appendix}

In this section we discuss chain rules for horizontal
derivatives that are used especially in connection with the morphism property.

\begin{prop}\label{p:chain} Let $\Heis$ be the standard sub-Riemannian Heisenberg group, and suppose that $(M,HM,g_M)$ is a contact sub-Riemannian manifold.
 Let $U$ be a domain in $\Heis$, $V$ a domain in $M$, and $f:U \to V$ a continuous function. Assume further that $u: V \to \mathbb{R}$ is smooth and $\Psi: V \to \Heis$ is a smooth chart so that $X(\Psi_i \circ f)$ and $Y(\Psi_i \circ f)$, $i\in \{1,2,3\}$, exist in the weak sense and belong to $L_{loc}^p$ for some $1\leq p<\infty$. Then the weak derivatives $X(u\circ f)$, $Y(u\circ f)$ exist, belong to $L_{loc}^p$ and are given almost everywhere by the following formulae
 \begin{align*}
  X(u\circ f)(p)=\sum_{i=1}^3\frac{\partial (u\circ \Psi^{-1})}{\partial \Psi_i}(\Psi(f(p))) X(\Psi_i \circ f)(p)
 \end{align*}
and
 \begin{align*}
  Y(u\circ f)(p)=\sum_{i=1}^3\frac{\partial (u\circ \Psi^{-1})}{\partial \Psi_i}(\Psi(f(p))) Y(\Psi_i \circ
  f)(p).
 \end{align*}
\end{prop}


\begin{proof}
We use charts to write locally $u\circ f= (u \circ \Psi^{-1})\circ (\Psi\circ f)$.
Note that $u\circ \Psi^{-1}$ is a smooth function on a domain in
$\Heis$ and thus the usual derivative $(u\circ
\Psi^{-1})_{\ast}$ exists everywhere in the domain of $u\circ
\Psi^{-1}$. Moreover, $u\circ \Psi^{-1}$ is locally Lipschitz both with respect to
the Euclidean and the sub-Riemannian metric on $\Heis$ and
thus $u\circ \Psi^{-1}$ is ACL (see Remark \ref{r:ACL_char}).

Concerning the factor $\Psi\circ f$, we denote
\begin{displaymath}
 \gamma_{\Psi \circ f,p}(s):= \Psi(f(p \exp(sX))).
\end{displaymath}
By assumption and Remark \ref{r:ACL_char}, for almost every $p$ in
the plane transversal to $X$, the tangent vector $
\dot{\gamma}_{\Psi \circ f,p}(s)$ exists for almost every $s$ and
it equals
\begin{align*}
X(\Psi \circ f)(p\exp(sX))= \sum_{i=1}^3 X(\Psi_i \circ
f)(p\exp(sX))
\partial_{\Psi_i}.
\end{align*}
Thus, for
\begin{displaymath}
 \gamma_{u\circ f,p}(s):= u(f(p \exp(sX)))= (u\circ \Psi^{-1})\circ (\Psi(f(p\exp (sX)))),
\end{displaymath}
we obtain
\begin{displaymath}
 X(u\circ f)(p \exp(sX))=\dot{\gamma}_{u\circ f,p}(s) =(u\circ \Psi^{-1})_{\ast,\Psi(f(p\exp(sX)))} X(\Psi \circ f)(p\exp(sX)),
\end{displaymath}
for almost every $s$, and analogously for $X$ replaced by $Y$.
This yields the formula for the chain rule. The fact that
$X(u\circ f)$ and $Y(u\circ f)$   belong to $L_{loc}^p$ is
immediate from the corresponding property of the horizontal
derivatives of $\Psi_i \circ f$ and the fact that $u\circ
\Psi^{-1}$ is smooth. Finally, we refer again to Remark
\ref{r:ACL_char} to deduce that the horizontal derivatives exist
also in a weak sense.
\end{proof}


Next we consider the case where the function $u$ is not smooth but only belongs to some Sobolev space. In this case we have to impose a stronger assumption on the map $f$, namely we will assume that it is quasiregular. For our purposes it suffices to discuss mappings between domains in the Heisenberg group.

\begin{prop}\label{p:chainQRHeis}
Let $f:\Omega \to \Omega'$ be a nonconstant quasiregular map between domains in $\Heis$, and let $u:\Omega' \to \mathbb{R}$ be an $HW_{loc}^{1,4}$-function. Then $u\circ f$ belongs to $HW_{loc}^{1,4}$. Moreover
\begin{displaymath}
\nabla_H (u\circ f)(p) = \left(D_H f(p)\right)^T \nabla_H u(f(p)),\quad\text{a.e. }p\in \Omega.
\end{displaymath}
\end{prop}

\begin{proof}
The proof goes along the same lines as in the Euclidean case, \cite[Theorem 14.28]{MR1207810}, using the fact that quasiregular mappings on the Heisenberg group are weakly contact and differentiable almost everywhere in the sense of Pansu \cite{Pan89}. Moreover, as shown in \cite[\S 5]{MR1778673}, the Pansu differential agrees almost everywhere with the map that is obtained by extending $D_H f$ to a homomorphism of the Lie algebra of $\Heis$.
\end{proof}

\bibliographystyle{acm}
\bibliography{references}
\end{document}